\documentclass[12pt,a4paper]{article}

\usepackage{a4wide,amssymb,amsmath,amsthm,xspace,ScriptFonts,epsfig}

\begin{document}

\newcommand{\RR}{\mathsf{R}}
\newcommand{\LL}{\mathsf{L}}

\def\co{\colon}

\newcommand{\R}{\Bbb R}
\newcommand{\Z}{\Bbb Z}
\newcommand{\Q}{\Bbb Q}
\newcommand{\C}{\Bbb C}
\newcommand{\esp}{\vskip .3cm \noindent}

\newcommand\Cfty{${\cal C}^{\infty}$}

\mathchardef\flat="115B

\newcommand\adm{{\subset}^a}
\newcommand\cfp{cofinal${}^+$}
\newcommand\cfm{cofinal${}^-$}
\newcommand\bdp{bounded${}^+$}
\newcommand\bdm{bounded${}^-$}

\def\CC#1{${\cal C}^{#1}$}
\def\h#1{\hat #1}
\def\t#1{\tilde #1}
\def\wt#1{\widetilde{#1}}
\def\wh#1{\widehat{#1}}
\def\wb#1{\overline{#1}}

\def\restrict#1{\bigr|_{#1}}

\newtheorem{thm}{Theorem}
\newtheorem*{thm*}{Theorem}
\newtheorem*{lemma*}{Lemma}
\newtheorem*{thmbis}{Theorem \ref{theorem1}'}

\newtheorem{thmalph}{Theorem}
\renewcommand\thethmalph{\Alph{thmalph}}

\newtheorem{ex}{Example}[section]

\newtheorem{prop}{Proposition}[section]
\newtheorem{lemma}[prop]{Lemma}
\newtheorem{cor}[prop]{Corollary}
\newtheorem{rem}[prop]{Remark}
\newtheorem{defi}[prop]{Definition}

\newtheorem{prob}{Problem}[section]
\newtheorem{q}{Question}[section]

\newtheorem*{prop*}{Proposition}
\newtheorem*{claim}{Claim}
\newtheorem*{lemma0}{Lemma 0}

\title{Homotopy in non metrizable $\omega$-bounded surfaces}
\author{Mathieu Baillif}
\date{\empty}

\maketitle

{\bf keywords:} {\em homotopy, non metrizable manifolds, longline}\\

{\bf Subject code classification:} {57N99, 14F35}

\abstract
{\footnotesize
We investigate the problem of describing the homotopy classes $[X,Y]$ of continuous functions between $\omega$-bounded
non metrizable manifolds $X,Y$. 
We define a family of surfaces $X$ built with the first octant $\mathsf{C}$
in $\LL^2$ ($\LL$ is the longline and $\RR$ the longray),
and show that $[X,\RR]$ is in bijection with so called `adapted' subsets of a partially ordered set. We also show
that $[M,\RR]$ can be computed for some surfaces $M$ that, unlike $\mathsf{C}$, do not contain $\RR$. 
This indicates that when 
$X,Y$ are $\omega$-bounded non metrizable surfaces, there might
be a link between $[X,Y]$ and the concept of $Y$-directions in $X$.
Many pictures are used and the proofs are quite
detailed.}
\endabstract

\section{Introduction}
Let $X,Y$ be surfaces. Can we describe the set $[X,Y]$ of homotopy classes of continuous maps $X\to Y$? 
If $X,Y$ are compact, the answer is known for a long time: a complete description of $[X,Y]$ is given by the morphisms
of the fundamental groups of $X,Y$ if $Y$ is not the sphere, and by the degree if $Y$ is the sphere (see for instance 
Theorem 11 page 428 in \cite{Spanier}). This paper investigates the case where $X,Y$ are $\omega$-bounded non metrizable
surfaces.\footnote{In fact, complete results will be given only when $Y$ is the longray, the simplest non metrizable
$1$-dimensional manifold.} The notion of $\omega$-boundedness is a kind of analogue to compactness in the non metrizable case.
This paper's aim is to describe in details some new phenomenas that appear when dropping the metrizability assumption,
and is a mix of original material and results that appeared in \cite{meszigues+Dave, meszigueshom}. The
proofs are elementary in the sense that apart from some knowledge about countable ordinals, no sophisticated theory is needed for most of the discussion (Section \ref{sec6} and Appendix \ref{sec5} are exceptions).

The class of $\omega$-bounded manifolds is interesting from a homotopical point of view because its
members cannot be both contractible and non metrizable (the proof will appear in \cite{meszigues+Nyikos}).
Therefore, if $X$ is a non metrizable $\omega$-bounded manifold, $|[X,X]|\ge 2$, for instance. (Notice that there are
contractible non metrizable surfaces, see \cite[Appendix A]{Spivak:vol1}.) Moreover, $\omega$-bounded {\em surfaces}
have been classified (in some way) by Nyikos in \cite{Nyikos:1984}: Each one consists of a compact metrizable `bag' 
($n$-torus with boundary) to which a
finite number of longpipes are attached. Longpipes will be defined below, they are `long' versions of the cylinder
$\mathbb{S}^1\times\R_{\ge 0}$. (Note: In contrast with the compact case, there are uncountably many non homeomorphic
longpipes.) 

This paper is example driven, we will not look for the greatest generality. Our aim is to 
show that in the class of $\omega$-bounded non metrizable manifolds, 
the above mentioned link between $[X,Y]$ and the $\pi_i(X),\pi_i(Y)$ may completely disappear:
The latter can all be trivial and the former quite complicated. 
In the other hand there seems to be a partial order structure behind $[X,Y]$  
(at least when $Y$ is the longray $\RR$). This partial order is defined on the set of homotopy classes of embeddings of $\RR$
in $X$. The examples we chosed are however all build with the same `brick' (Section \ref{sec6} is again 
an exception), the {\em first octant} $\mathsf{C}$ to be defined below, which contains two homotopy classes of embeddings
of $\RR$, and this partial order structure is perhaps a feature of this brick. Meanwhile, the study (in Section \ref{sec6})
of a completely
different example containing {\em no} copy of $\RR$ indicates that there might be a general phenomenon linked with the
concept of {\em direction} (see the end of Section \ref{sec5}). 
Homotopy in non metrizable manifolds is a relatively 
new subject: a lot of things remain to be done, and we hope this paper will contribute to its popularization.
We indicate along the way some among the many problems we thought worth of further study.

The paper is organized as follows. Basic facts such as the definition of the longray, of $\mathsf{C}$, etc., are given
in Section \ref{Basics}. In Section \ref{sec2} we define a family of $\omega$-bounded manifolds built with
(countably many) copies of $\mathsf{C}$ and the partial order associated, and we state Theorem \ref{theorem1} linking
$[X,\RR]$ to this partial order. Section \ref{sec3} contains the `technical lemmas' about $\mathsf{C}$, mainly partition
properties of continuous maps $\mathsf{C}\to\RR$. These properties are the key point of the proof of Theorem \ref{theorem1},
a proof we complete in Section \ref{sec4}. Section \ref{sec6} deals with quite different manifolds obtained from the tangent
bundle of $\RR$. We shall use some results of Nyikos \cite{Nyikos:1992}. Finally, we show in Appendix \ref{sec5} how
the construction of Section \ref{sec3} could be pushed further to obtain manifolds with uncountably many homotopy
classes of embedding of $\RR$, and how to obtain similar results. We shall provide much less details there.

Notice that there are two classical ways of defining an homotopy between two maps
$f,g\co X\to Y$: The first is a continuous $\phi\co X\times[0,1]\to Y$ with $\phi(\cdot,0)=f$, $\phi(\cdot,1)=g$, and
the second is a continuous $\varphi\co [0,1]\to\mathcal{C}(X,Y)$ with $\varphi(0)=f$, $\varphi(1)=g$, where
$\mathcal{C}(X,Y)$ is the space of continuous maps $X\to Y$ with the compact-open topology.
These two notions are equivalent if $X$ is locally compact, a proof can be
found in the appendix of \cite{Hatcher}.


\section{Basics}\label{Basics}
We follow von Neumann's definition of ordinals, identifying $\alpha$ with the set of ordinals
smaller than it. Thus, the first uncountable ordinal $\omega_1$ is the set of (finite or) countable ordinals, and 
$\omega$ the set of finite ordinals, i.e. the natural numbers.
Recall that if $\{\alpha_m\}_{m\in\omega}$ is a countable set of countable ordinals, then
$\sup_{m\in\omega}\alpha_m$ is also a countable ordinal. $\R$ denotes the real numbers.

By a {\em manifold} is meant a connected Hausdorff space, each of whose points possess a neighborhood homeomorphic to $\R^n$ for
some $n$ (which is fixed, by connectedness, and is called the dimension of the manifold). A {\em surface} is a $2$-dimensional
manifold. 
We allow manifolds with boundary (amending the definition consequently).
Recall that in the category of manifolds, metrizability is equivalent to
seemingly weaker properties as Lindel\"ofness, second countability, paracompactness, etc. (see 
\cite{Gauld:Met} for an impressive list). A topological space $X$ is called {\em type I} if  
\begin{equation}
  \label{type I}
  X=\bigcup_{\alpha<\omega_1}U_\alpha,
\end{equation}
where $U_\alpha$ is open, $\wb{U}_\alpha$ Lindel\"of, $\wb{U}_\alpha\subset U_\beta$ when $\alpha<\beta$.
A type I manifold $X$ is Lindel\"of (and thus metrizable) if and only if $X=U_\alpha$ for some $\alpha<\omega_1$.
Thus, the non metrizability of a type I manifold comes from its `wideness' rather than from its `shape'. A manifold is
$\omega$-bounded if it is type I and sequentially compact (or equivalently if the closure of any countable set is compact,
which is the `official' definition, see  
Corollary 5.4 in \cite{Nyikos:1984}). 

We recall that the (closed) long ray is $\RR=\omega_1\times[0,1[$ with the lexicographic order $\le$
and the order 
topology. In other words, we glue together $\omega_1$ copies of $[0,1[$. 
(Notice that $\R_{\ge 0}$ is homeomorphic to $\omega\times[0,1[$ with the lexicographic order topology.)
To simplify notation, we shall
denote $(\alpha,0)\in\RR$ simply by $\alpha$, and often treat $\omega_1$ as a subset of $\RR$.
The open long ray is $\RR\backslash\{0\}$. We define $[x,y]$, $]x,y[$ (and so on) in the usual way on any totally ordered set 
($\RR$ and $\omega_1$ for instance). 
To see that $\RR$ is a $1$-dimensional manifold (with boundary) 
one considers the atlas $U_\alpha=[0,\alpha[\subset\RR$, and show by induction that
$U_\alpha$ is homeomorphic to $[0,1[\subset\R$, using the fact that for any countable limit ordinal
$\alpha\in\omega_1$, there is a sequence $\alpha_m<\alpha$ ($m\in\omega$) converging to $\alpha$.
This also shows that $\RR$ is of type I.
With a little more effort one sees that $\RR$ can be given a structure of \CC{\infty} or even of analytic manifold,
see \cite{Nyikos:1992} and references therein.
Since $\RR$ is non Lindel\"of (for instance, the atlas gives a cover with no countable subcover),
$\RR$ is non metrizable, non paracompact, and so on; it is moreover non separable. 
However, $\RR$ is sequentially compact (and thus $\omega$-bounded): Any sequence $x_m\in\RR$ ($m\in\omega$)
is contained in $\wb{U}_\alpha$
for $\alpha=\sup_{m\in\omega}\alpha_m$, where $\alpha_m$ is such that $x_m\in U_{\alpha_m}$;
and $\wb{U}_\alpha$ is homeomorphic to $[0,1]$ and thus sequentially compact.
We shall use sequential compactness of $\RR$ (and of surfaces built with $\RR$) thoroughly throughout this
paper when saying that some (sub)sequence converges, most of the time without mentioning it.

It is a good exercise to show that $\RR$ is non contractible, due to its `wideness'; however, each
$\pi_i(\RR)$ is evidently trivial (the continuous image of a compact set must be contained in some 
$U_\alpha\simeq [0,1[$), as will be those of all the surfaces we shall consider. 
\begin{defi}\label{defclub}
  We say that a subset of $\RR$ or $\omega_1$ is club if it is closed and unbounded.
\end{defi}
The following lemma, whose proof can be found for instance in \cite{Kunen},
shows the importance of club subsets:
\begin{lemma}\label{club}
  If $E_m\subset\RR,\omega_1$ are club for each $m\in\omega$, so is $\bigcap_{m\in\omega}E_m$.
\end{lemma}
The word `cofinal' will be used as a synonym of `unbounded' for maps $\RR\to\RR$.
The following lemma is, so to say, the `canonical representative' of many results that we shall give below;
it says that when $f_k\co \RR\to\RR$ ($k\in\omega$) are cofinal and continuous, then each
$f_k$ preserves some `blocks'
$[\beta_\gamma,\beta_{\gamma+1}]$. So, if $f,g\co \RR\to\RR$ are such functions, it is very easy
to prove that $f$ and $g$ are homotopic: We just define homotopies relative to the boundary in
each $[\beta_\gamma,\beta_{\gamma+1}]$ (that is, we apply Lemma \ref{homotopy} below).\footnote{This gives
a proof that $[\RR,\RR]\simeq\{0,1\}$, a fact already proved by D. Gauld in \cite{Gauld}. } 
This reduction to a trivial case is very typical of the methods we shall use throughout this paper.
Because similar ideas will be used over and over, we chosed to include next lemma's proof, eventhough
it already appeared in \cite{meszigues+Dave}.
We shall
abusively write $[x,\omega_1[$ for the set of $y\in\RR$ that are $\ge x$.
\begin{lemma}
\label{lemme1} Let $f_k\co \RR\to\RR$ be continuous for all $k\in\omega$.
If each $f_k$ is bounded, there is an $x\in\RR$ such that $f_k\restrict{[x,\omega_1[}$ is constant $\forall k$.
If $f_k$ is cofinal for all $k$, there is a $\omega_1$-sequence $\beta_\gamma\in\omega_1\subset\RR$ with 
$\gamma<\gamma'\Rightarrow \beta_\gamma<\beta_{\gamma'}$ and 
$\beta_\gamma=\lim_{\gamma'<\gamma}\beta_{\gamma'}$ if $\gamma$ is a limit ordinal, such that
$f_k([\beta_\gamma,\beta_{\gamma+1}])=[\beta_\gamma,\beta_{\gamma+1}]$ if $\gamma>0$,
and $f_k([0,\beta_0])\subset [0,\beta_0]$.
\end{lemma}
\proof
Assume that each $f_k$ is bounded.
It is easy to check (by hand or applying \cite[Theorem 3.10.6]{Engelking})
that $f_k$ attains its bounds.
Define now an increasing sequence $y^k_m$ such that $y^k_0=0$, $y^k_{2m+1}$ gives the supremum of
$f_k\restrict{[y^k_{2m},\omega_1[}$ and $y^k_{2m}$ the infimum of
$f_k\restrict{[y^k_{2m-1},\omega_1[}$. 
Letting $x_k=\sup_{m\in\omega}y^k_m$ implies that
$f_k\restrict{[x_k,\omega_1[}$ is constant. (By the way, this is \cite[Exercise (42), p. 91]{Kunen}.)
Take $x=\sup_{k\in\omega} x_k$.
If $f_k$ is cofinal, for all $\alpha$ there is a $\gamma(\alpha)\ge\alpha$ such that 
$f_k([\gamma(\alpha),\omega_1[)\subset[\alpha,\omega_1[$. (Otherwise the set of 
$\gamma\ge\alpha$ such that $f_k(\gamma)\le\alpha$ is club, but so is the set of 
$\gamma'\ge\alpha$ with $f_k(\gamma')\ge\alpha+1$, which is a contradiction by Lemma \ref{club}.)
It follows that $E_k=\{\gamma\in\omega_1\subset\RR\,:\,f_k([\gamma,\omega_1[)\subset[\gamma,\omega_1[\}$ is
club: closeness is immediate, for unboundedness, start with $\alpha_0$ and define
$\alpha_{m+1}=\gamma(\alpha_m)$; its limit is in $E$. 
But $F_k=\{\gamma\in\omega_1\subset\RR\,:\,f_k([0,\gamma])\subset[0,\gamma]\}$ is also club: For unboundedness,
given some $y_0^k$, define $y_m^k=\max\{y_{m-1}^k,\sup_{[0,y^k_{m-1}]}f_k\}$ ($m\in\omega$), then $y_m^k\to x_k\in F_k$ by
continuity.
Hence, $D=\bigcap_{k\in\omega}E_k\cap F_k$ is club. We define $\beta_\gamma$ by transfinite induction, choosing
$\beta_0=\min D$, $\beta_{\gamma+1}=\min D\cap [\beta_\gamma+1,\omega_1[$, and if $\gamma$ is limit,
$\beta_{\gamma}=\sup_{\gamma'<\gamma}\beta_{\gamma'}$.
\endproof

Notice that $[0,\beta_0]\cup\bigcup_{\gamma\in\omega_1}[\beta_\gamma,\beta_{\gamma+1}]=\RR$.
We end this section with
the only purely homotopical lemma that we shall use in this paper, whose proof is trivial:
\begin{lemma}\label{homotopy}
  Let $f,g\co X\to Y$ be continuous and $Y$ be homeomorphic to $[0,1]^d$. Then,
  there is a homotopy $h_t$ such that $h_0=f,h_1=g$ and for all $t$,
  $h_t\bigr|_Q=id$, where $Q=\{x\in X\,:\,f(x)=g(x)\}$.
\end{lemma}


\section{Some (simple) $\omega$-bounded non metrizable surfaces}\label{sec2}
A space is a {\em longplane} if it is the union of a chain (\ref{type I})
where $U_\alpha$ is open and homeomorphic to $\R^2$,
$\wb{U}_\alpha\subset U_\beta$ and the boundary of $U_\alpha$ in $U_\beta$ is homeomorphic to the circle $\mathbb{S}^1$
when $\alpha<\beta$. A longpipe is a longplane with a point removed (see \cite{Nyikos:1984}). We prefer to work with 
longplanes because they have trivial homotopy groups.
We now define the `building brick' of almost all surfaces of this paper:
\begin{defi}\label{defC}
  We let the first octant be 
  $\mathsf{C}=\{(x,y)\in\RR^2\,:\,y\le x\}$, with the induced topology.
\end{defi}
As said before, the homotopy properties of the manifolds that we will define with $\mathsf{C}$ are closely related
to some embeddings (or copies) of $\RR$ in the manifolds.
In $\mathsf{C}$, there are (many) horizontal and (one) diagonal `canonical' copies of $\RR$:
\begin{defi}
  Let $\Delta_d=\{(x,x)\,:\, x\in\RR\}$.
  For $c\in\RR$, we let $\Delta_h(c)$ be $\{(x,c)\,:\, x\in\RR, x\ge c\}$, and
  write $\Delta_h(0)$ as $\Delta_h$.
\end{defi}
(Of course, $d$ and $h$ stand for `diagonal' and `horizontal'.)
One sees immediately that $\Delta_h(c)$ and $\Delta_d$ with the induced topology are homeomorphic to $\RR$.
The two following lemmas, whose proofs are left as good exercises (or can be found in \cite{Nyikos:1984}),
show that $\Delta_d,\Delta_h(c)$ are `topologically different'.
\begin{lemma}\label{open sets}
  Let $U\subset\RR^2$ be open.
  If $\Delta_d\subset U$, there is an $x\in\RR$ such that $[x,\omega_1[^2\subset U$.
  If $\Delta_h(c)\subset U$, there are $y<c<y'$ and $x$ in $\RR$ with $U\supset [x,\omega_1[\times]y,y'[$.
\end{lemma}

\begin{lemma}\label{DeltaC}
  If $\Delta\subset\mathsf{C}$ is a copy of $\RR$, then either there is a $c$ for which 
  $\Delta$ is contained in $\Delta_h(c)$ 
  outside of a compact set, or $\Delta\cap\Delta_d$ is club.
\end{lemma}
\ \\
\parbox{0.6\textwidth}{When we picture $\mathsf{C}$, we use an arrow pointing at $\Delta_d$ (see
the figure opposite) 
as a graphic tool that makes
explicit which boundary is $\Delta_h$ and which is $\Delta_d$. This graphic convention has also an
homotopic meaning (Lemma \ref{h d cof}).}
\parbox{0.05\textwidth}{\ }
\parbox{0.3\textwidth}{\epsfig{figure=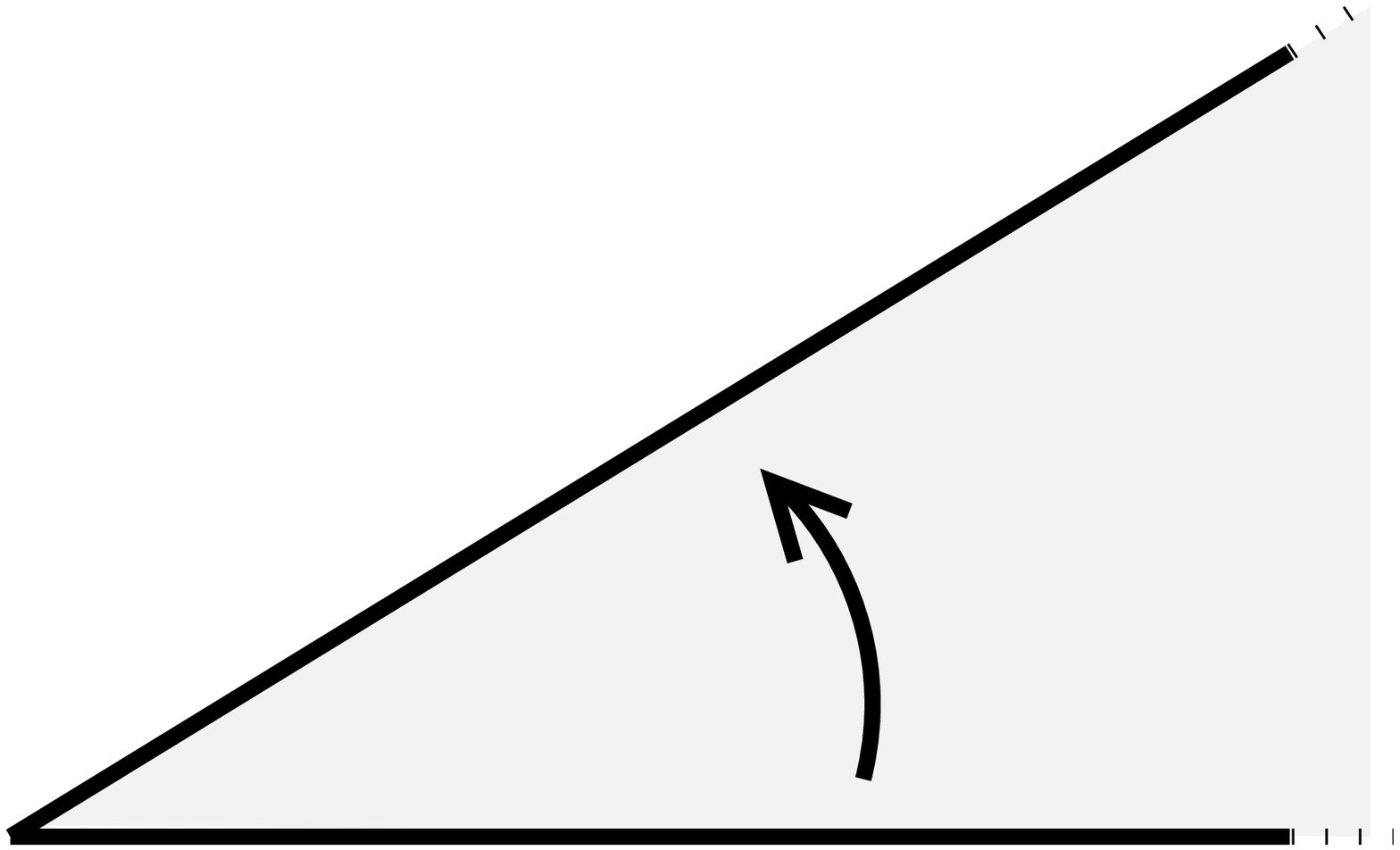,height=2cm}}
There are several ways to glue together two copies $C_0,C_1$
of $\mathsf{C}$ along their boundary components.
We will be interested in those where $C_1$ is glued `on the top of $C_0$',
as on the picture below.
\begin{center}
  \epsfig{figure=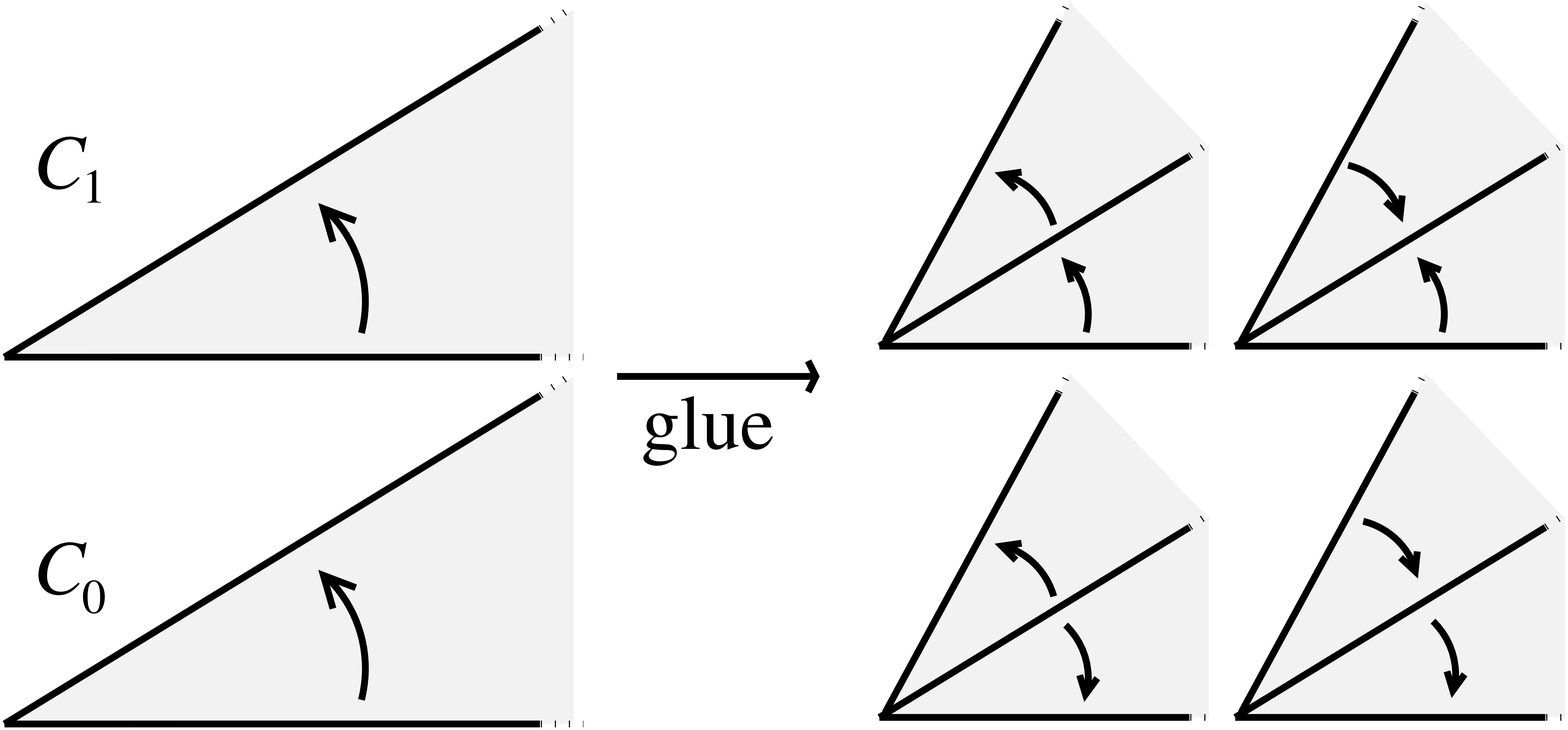,height=4cm}
\end{center}
We symbolize these gluings by the following pairs of $\uparrow,\downarrow$: 
$\left<\uparrow\uparrow\right>$, $\left<\uparrow\downarrow\right>$,
$\left<\downarrow\uparrow\right>$, and $\left<\downarrow\downarrow\right>$.
Given a finite sequence $s\co \{0,\dots,n-1\}\to\{\uparrow,\downarrow\}$, we define the surface $M_{n,s}$ 
by induction, gluing copies $C_i$ ($i=0,\dots,n-1$) of $\mathsf{C}$ `on the top of each other',
coherently with the sequence $s$. Then, we glue together the boundary components (copies of $\RR$) that
remained free. 
We denote the copies of $\Delta_h,\Delta_d\subset C_i\subset M_{n,s}$ by
$\Delta_i$ ($i=1,\dots,n$),
turning counterclockwise, as on the examples below.
\begin{center}
  \epsfig{figure=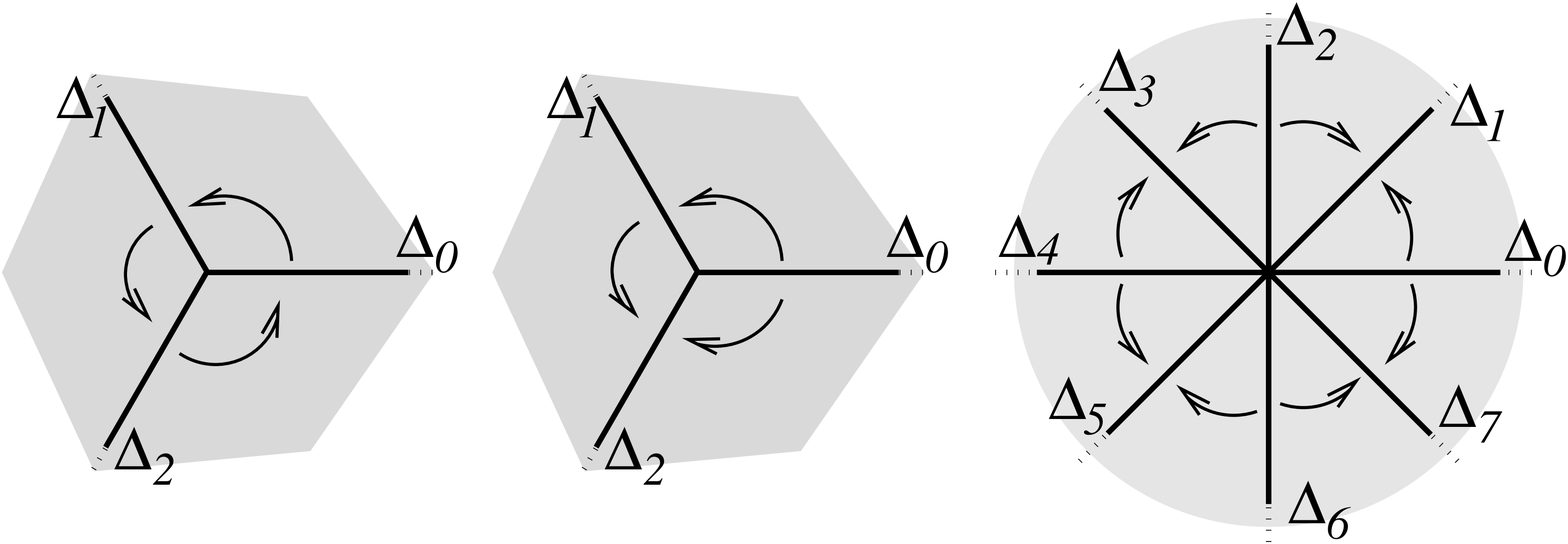,height=3.5cm} \\
  {\em Figure 1: The surfaces $M_{3,\left<\uparrow\uparrow\uparrow\right>}$, 
  $M_{3,\left<\uparrow\uparrow\downarrow\right>}$ and 
  $M_{8,\left<\uparrow\downarrow\uparrow\downarrow\uparrow\downarrow\uparrow\downarrow\right>}$.}
\end{center}
(Notice that the rightmost example is $\LL^2$, where $\LL$ is the longline, which consists in two copies of
$\RR$ glued at $0$.) \\
\parbox{0.5\textwidth}{
Now, if $s\co \omega\to\{\uparrow,\downarrow\}$ is an {\em infinite} sequence, we may define $M_{\omega,s}$ 
in the same way,
with $\omega$ copies $C_i$ of $\mathsf{C}$,
with the $C_i$ accumulating on $\Delta_0$ (that is, we give $x\in\Delta_0$ a neighborhood basis consisting of 
unions $U_0\cup U_m\cup \bigcup_{i\ge m+1}U_i$, where
the $U_i\subset C_i$ are 
as pictured opposite). Of course, $s$ determines
in which way each $C_i$ is glued. }
\parbox{0.45\textwidth}{\epsfig{figure=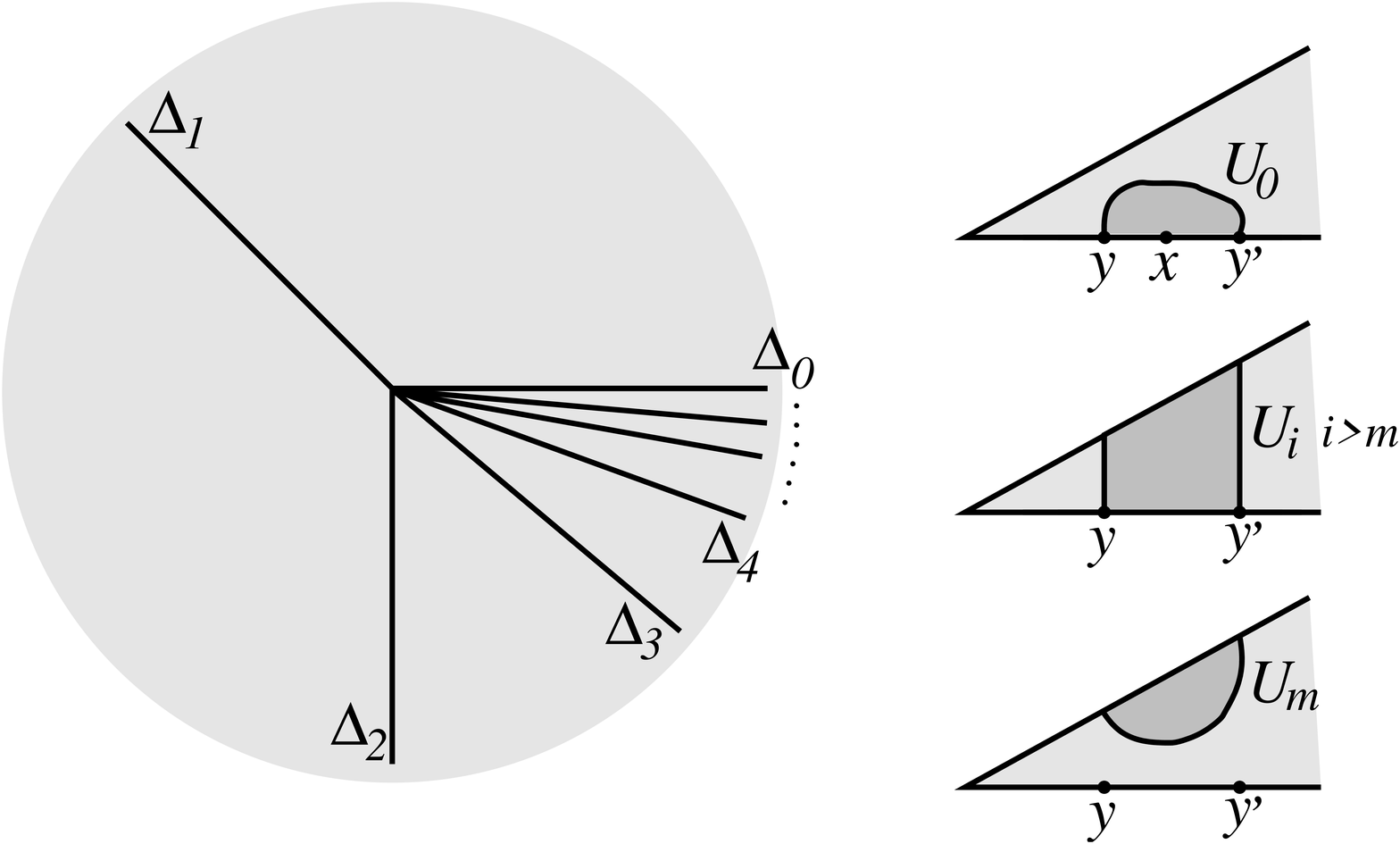,height=3.5cm}}
\\
There is no reason to stop at $\omega$, and we define similarly
$M_{\alpha,s}$ for any $s\co \alpha\to\{\uparrow,\downarrow\}$,
with $\alpha$ any limit ordinal $<\omega_1$. (If $\alpha=\beta+n$ with $n\in\omega$, $\beta$ limit,
and $s\co \alpha\to\{\uparrow,\downarrow\}$, 
then $M_{\alpha,s}$ is homeomorphic to $M_{\beta,s'}$, where $s'\co \beta\to\{\uparrow,\downarrow\}$,
$s'(i)=s(\beta+i)$ if $i\le n$, $s'(n+i+1)=s(i)$, and $s'(\gamma)=s(\gamma)$ for $\omega\le\gamma<\beta$. 
We may thus forget about successor ordinals.)

It should be clear that all the $M_{\alpha,s}$ we defined are longplanes, and in particular,
$\omega$-bounded surfaces 
with trivial homotopy groups.

We have defined $M_{\alpha,s}$ for $s\co \alpha\to\{\uparrow,\downarrow\}$, $\alpha<\omega_1$. We 
now define a partial ordering on $\alpha$ associated with $s$.
\begin{defi}\label{po} 
  Let $s\co \alpha\to\{\uparrow,\downarrow\}$, $\alpha<\omega_1$.
  We let $\mathbb{P}_{\alpha,s}$ be the partially ordered set $\left<\alpha\,:\,\prec\right>$ where $\prec$ 
  is the reflexive and
  transitive closure
  of $\prec'$ defined by 
  $$
    \begin{array}{lcr}\gamma\prec'\gamma+1&\text{ if }&s(\gamma)=\uparrow, \\ 
      \gamma+1\prec'\gamma&\text{ if }&s(\gamma)=\downarrow,\end{array}
  $$
  and, if $\alpha=n\in\omega$,
  $$
    \begin{array}{lcr}n-1\prec' 0&\text{ if }&s(n-1)=\uparrow, \\ 
      0\prec'n-1&\text{ if }&s(n-1)=\downarrow.\end{array}
  $$
\end{defi}
This partial order can be seen on pictures: $\gamma\prec\gamma'$ if we can go from $\Delta_\gamma$
to $\Delta_{\gamma'}$ following a finite number of arrows. Notice that $\prec$ may not be a partial order
in the strict sense: If $\alpha=n\in\omega$ and $s(i)=\uparrow$ for all $i$, then $i\prec j\prec i$ for all $i,j$.
For instance, the two leftmost examples of Figure 1 give the orders
$0\prec 1\prec 2\prec 0$ for $s=\left<\uparrow\uparrow\uparrow\right>$ and
$0\prec 1\prec 2\not\prec 0$, $1\not\prec 0$, for $s=\left<\uparrow\uparrow\downarrow\right>$.
\begin{defi}\label{adapted}
  Let $s\co \alpha\to\{\uparrow\,\downarrow\}$, $\alpha\le\omega_1$. We say that $W\subset\mathbb{P}_{\alpha,s}$ is adapted
  if 
  \begin{itemize}
  \item[i)] $\gamma\in W$, $\gamma\prec\gamma'$ $\Rightarrow$ $\gamma'\in W$,
  \item[ii)] $\forall\beta\le\alpha$ with $\beta$ limit,
  $\exists\gamma(\beta)<\beta$ such that
  $$
    \begin{array}{l} 0\in W \text{ if }\alpha=\beta \\
    \beta\in W\text{ if }\beta<\alpha\end{array}\quad\Leftrightarrow\quad
    \gamma(\beta)\in W \quad\Leftrightarrow\quad \gamma'\in W\,\forall\gamma'\in [\gamma(\beta),\beta[.
  $$
  \end{itemize}
\end{defi}
Notice that condition ii) is empty if $\alpha<\omega$. The theorem that links $\mathbb{P}_{\alpha,s}$ with
$[M_{\alpha,s},\RR]$ is the following:
\begin{thm}\label{theorem1}
  Let $\alpha <\omega_1$, $s\co \alpha\to\{\uparrow\,\downarrow\}$, and $M_{\alpha,s}$, 
  $\mathbb{P}_{\alpha,s}$ be as above. Then,
  $[M_{\alpha,s},\RR]\simeq\{ W\subset\mathbb{P}_{\alpha,s}\,:\,W\text{ adapted}\}$.
\end{thm}
Here $\simeq$ means that there is a natural bijection between the two sets.
If $\alpha<\omega$, we simply recover the following result of \cite{meszigueshom}:
\begin{thm*}
  Let $s:n\to\{\uparrow\,\downarrow\}$, then 
  $[M_{n,s},\RR]\simeq\{ W\subset\mathbb{P}_{n,s}\,:\,W\text{ is an antichain}\}$.
\end{thm*}  
(Recall that an antichain is a set of incomparable elements. Given $W$ adapted, its minimal elements form 
an antichain.)
Theorem \ref{theorem1} will be proved in the next two sections.
\begin{rem}
  It is perhaps intructive to notice that there is a rough analogy between (compact) metrizable surfaces and 
  those considered in this section (and in Appendix \ref{sec5}). There is a huge amount of theory involving homology
  and homotopy in the class of CW-complexes (see for instance \cite{Spanier}), and
  one of its interests lies in the a priori non evident fact that
  any compact manifold has the homotopy type of a CW-complex (and can be triangulated in dimension $<4$).
  Here, $\mathsf{C}$ can be thought as a kind of non metrizable analogue of the $2$-dimensional simplex, with which we
  have built our surfaces. But a lot of things are different, though: $\mathsf{C}$ is non contractible, and more importantly
  many surfaces (those of Section \ref{sec6} for instance, or more simply the long cylinder $\mathbb{S}^1\times\RR$)
  cannot be ``$\mathsf{C}$-ulated''.
\end{rem}


\section{A close investigation of $\mathsf{C}$}\label{sec3}
(The results of this section appeared (in a more general form) in \cite{meszigueshom}, except
Propositions \ref{prop1}--\ref{prop2} when $\alpha\ge\omega$.) 
If $f\co \RR\to\RR$ is cofinal, it is not homotopic to a bounded function.
Thus, since there are homotopies sending $\Delta_h(c)$ to $\Delta_h(c')$, $f\restrict{\Delta_h(c)}$ is
unbounded if and only if $f\restrict{\Delta_h(c')}$ is. This motivates the following definition:
\begin{defi}
  For $*=d,h$, we say that $f\co \mathsf{C}\to\RR$ is 
  $*$-cofinal (resp. $*$-bounded) if 
  $f\restrict{\Delta_*}$ is unbounded (resp. bounded). 
\end{defi}
\begin{lemma}\label{h d cof}
  Let $f\co \mathsf{C}\to\RR$ be continuous. If $f$ is $h$-cofinal, it is $d$-cofinal.
\end{lemma}
\noindent
\parbox{0.6\textwidth}{\proof Let $z\in\RR$, we shall find some $x\in\Delta_d$ with $f(x)\ge z$.
Since $f\restrict{\Delta_h(0)}$ is cofinal, 
$f\restrict{\Delta_h(c)}$ is also cofinal for all $c$.
We define $x_m=(x^1_m,x^2_m)\in\mathsf{C}$, $m\in\omega$ inductively as follows. Let
$x_0\in\Delta_h(0)$ be such that $f(x_0)\ge z$, and choose 
$x_{m+1}\in\Delta_h(x^1_m)$ with $f(x_m)\ge z$. Then, (a subsequence of) $x_m$ converges to
some $x\in\Delta_d$ with $f(x)\ge z$, see opposite.\endproof}
\parbox{0.05\textwidth}{\ }
\parbox{0.35\textwidth}{\epsfig{figure=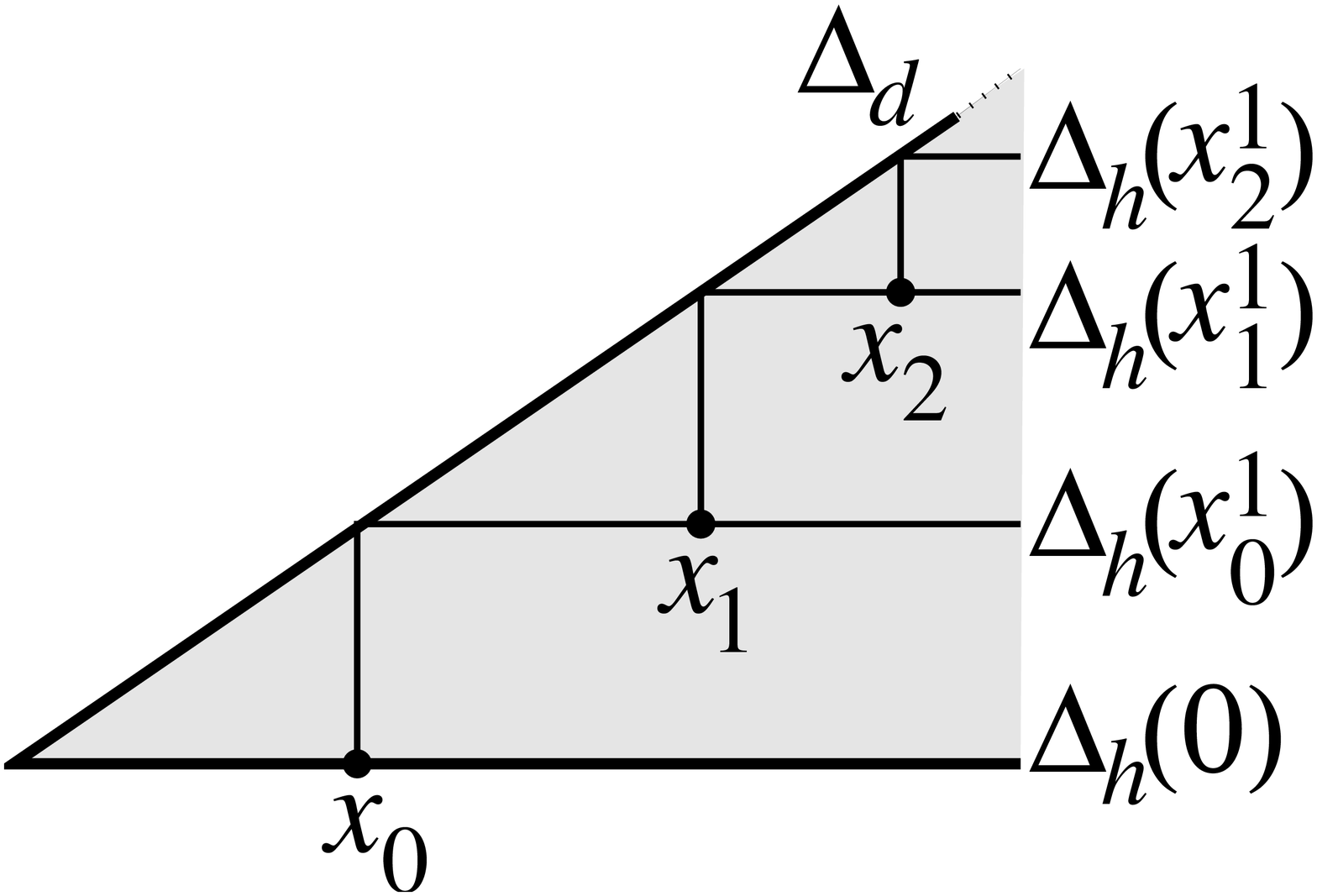,height=2.5cm}}

\begin{lemma}\label{d cofinal}
  Let $f\co \mathsf{C}\to\RR$ be continuous and unbounded. Then, $f$ is $d$-cofinal.
\end{lemma}
\proof
Suppose that $f\restrict{\Delta_d}$ is strictly bounded by $b$. Then
$\Delta_d\subset f^{-1}([0,b[)=U$, which is open. By Lemma \ref{open sets}, there is some
$x$ such that $f([x,\omega_1[^2\,\cap\,\mathsf{C})\subset [0,b]$. Let $E=(\RR\times[0,x])\cap\mathsf{C}$.
Since $f$ is unbounded, $f\restrict{E}$ is unbounded. 
Let $\{e_m\}_{m\in\omega}$ be a dense subset of $[0,x]$.
If $f\restrict{(\RR\times\{e_m\})\cap\mathsf{C}}$ is bounded by $b_m$, $f\restrict{E}$ is bounded by $\sup_m b_m$.
There is thus some $m$ such that $f\restrict{(\RR\times\{e_m\})\cap\mathsf{C}}$ is unbounded,
so $f$ is $h$-cofinal, and therefore $d$-cofinal
by Lemma \ref{h d cof}.
\endproof
There are thus three possible classes of cofinality for continuous functions $\mathsf{C}\to\RR$:
bounded functions, $h$-cofinal (and thus $d$-cofinal) 
functions (for example, the horizontal projection) and $d$-cofinal and $h$-bounded
functions (for example, the vertical projection).  
\begin{defi} Let $s\co \alpha\to\{\uparrow,\downarrow\}$.
  We say that a continuous $f\co M_{\alpha,s}\to\RR$ is $\gamma$-cofinal (resp. $\gamma$-bounded)
  if $f\restrict{\Delta_\gamma}$ is cofinal (resp. bounded).
  The cofinality class of $f$ is then
  $\mathfrak{C}(f)=\{\gamma<\alpha\,:\,f\text{ is }\gamma\text{-cofinal.}\}$.
\end{defi}
Half of Theorem \ref{theorem1} is provided by:
\begin{prop}\label{prop1} Let $s\co \alpha\to\{\uparrow,\downarrow\}$. Then,
  $$\{\mathfrak{C}(f)\,\,:\,\,f\co M_{\alpha,s}\to\RR\text{ continuous}\}=
  \{ W\subset\mathbb{P}_{\alpha,s}\,:\,W\text{ adapted}\}.$$
\end{prop}
Proposition \ref{prop1} explains point ii) of Definition \ref{adapted}: if $\beta$ is limit,
the $\Delta_\gamma$ for $\gamma<\beta$ accumulate on $\Delta_\beta$,
so $f$ is $\beta$-cofinal if and only if there is some
$\gamma(\beta)<\beta$ with $f$ $\gamma'$-cofinal for all $\gamma\le\gamma'<\beta$. 
Point i) comes from Lemma \ref{h d cof}.
\proof[Proof of Proposition \ref{prop1}]
  Let $s\co \alpha\to\{\uparrow,\downarrow\}$, $\alpha <\omega_1$ and
  $f\co M_{\alpha,s}\to\RR$ be continuous. By Lemma \ref{h d cof} and the definition of $\prec$, if
  $\gamma\in\mathfrak{C}(f)$ and $\gamma\prec\gamma'$, then $\gamma'\in\mathfrak{C}(f)$. 
  Together with the above remarks, this shows that 
  $\mathfrak{C}(f)$ is adapted. Conversely, given an adapted $W$, we may find
  $f$ with $\mathfrak{C}(f)=W$ by choosing $f\restrict{C_\gamma}$ 
  to be either $\equiv 0$, the vertical or the horizontal projection, according to $W$. 
  Since $W$ is adapted, $f$ will be continuous.
\endproof
To obtain Theorem \ref{theorem1}, it is therefore enough to prove the following:
\begin{prop}\label{prop2}
  Let $s\co \alpha\to\{\uparrow,\downarrow\}$, and
  $f,g\co M_{\alpha,s}\to\RR$ be continuous with $\mathfrak{C}(f)=\mathfrak{C}(g)$. 
  Then, $f$ and $g$ are homotopic.
\end{prop}
This will be done in the next section, using partition properties (`preservation of blocks') that we
now explain. Our goal is to find a closed $\omega_1$-sequence $\{\beta_\gamma\,:\,\gamma\in\omega_1\}$
which satisfies what is described in Figure 2 below.

\begin{center}
  \epsfig{figure=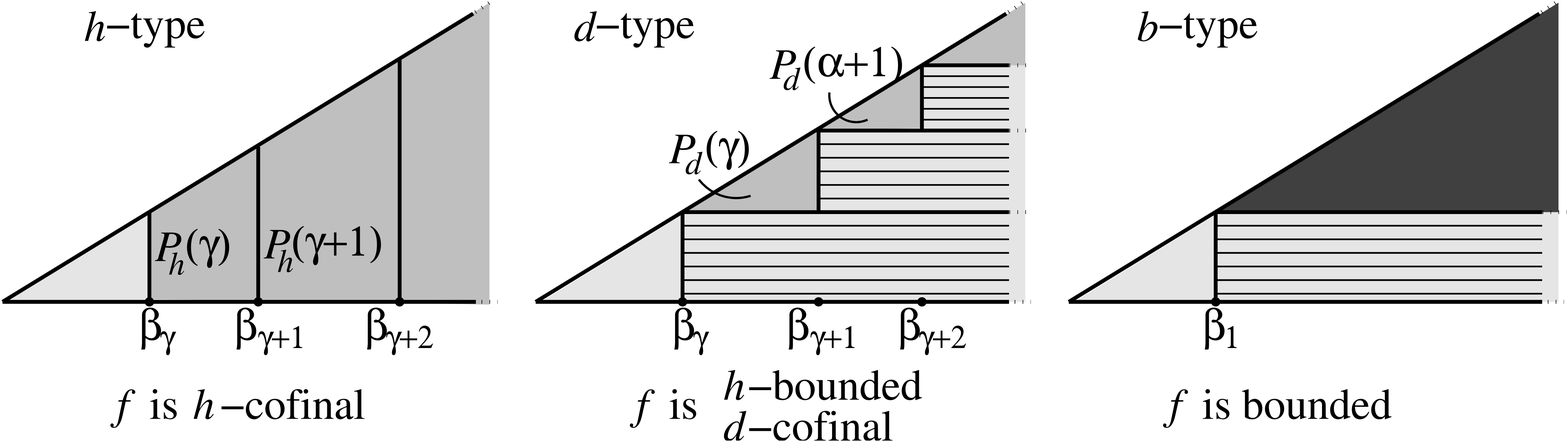,width=14cm}
  \\
  {\em Figure 2: Partition properties.}
\end{center}

The picture should be understood as follows. First, the partition for $f$ depends on its cofinality type, as
indicated under each figure. The grey regions denoted by $P_d(\gamma)$ or $P_h(\gamma)$ are sent 
by $f$ in $[\beta_\gamma,\beta_{\gamma+1}]$, and $f$ is constant on the horizontal thin lines and on the
black triangle at the righthandside. These properties will easily yield Proposition \ref{prop2} as we shall
see in the next section.
To prove them, we will show that some subsets of $\omega_1$ linked with the partition are club.
First, a definition of the `blocks'.
\begin{defi}
  Let $\gamma\in\omega_1\subset\RR$. We define
  $A^+_h(\gamma)=\{(x,y)\in\mathsf{C}\,:\, x\ge\gamma\}$,
  $A^+_d(\gamma)=[\gamma,\omega_1[^2\cap\mathsf{C}$,
  $A^-_d(\gamma)=A^-_h(\gamma)=\{(x,y)\in\mathsf{C}\,:\, x\le\gamma\}$ (see below).
  \begin{center}{\epsfig{figure=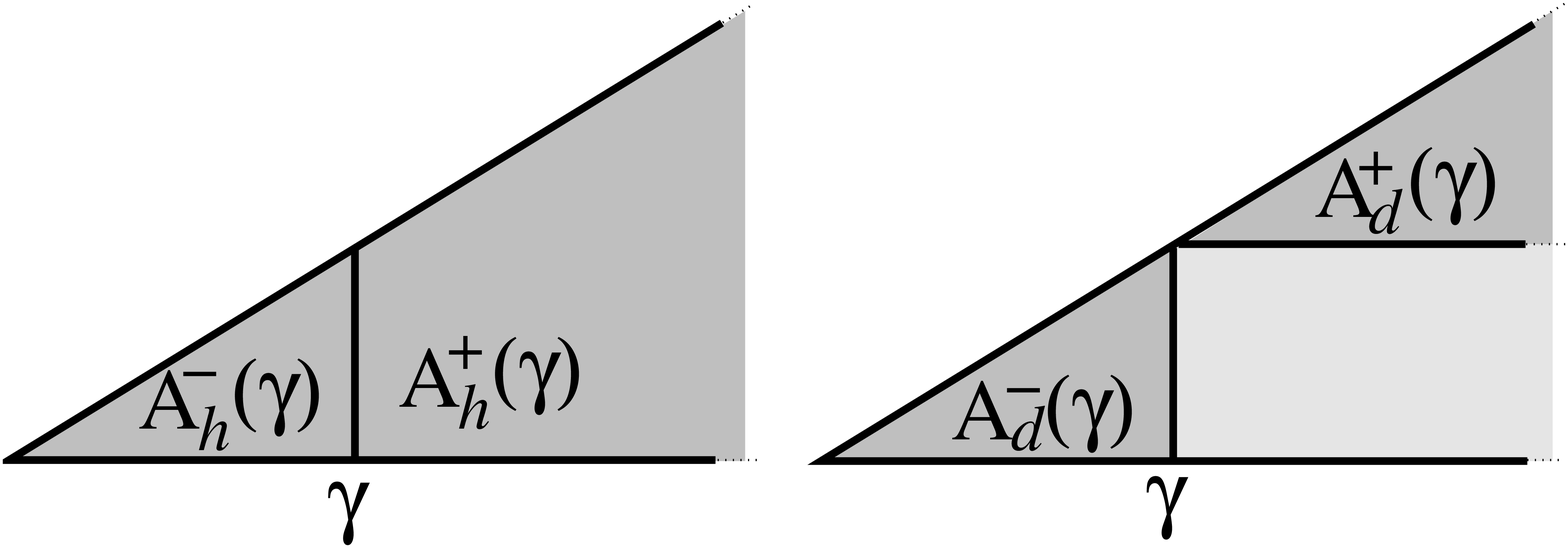,height=2.5cm}}\end{center}
\end{defi}
(Notice that in Figure 2, $P_*(\gamma)=A^+_*(\beta_{\gamma+1})\cap A^-_*(\beta_\gamma)$ for $*=d,h$.)
\begin{lemma}\label{equiv lemme 5.4}
  Suppose that $f\co \mathsf{C}\to\RR$ is continuous and $h$-cofinal. Then, for all 
  $c\in\RR$, there is $d(c)$ minimal such that 
  $f([d(c),\omega_1[\times[0,c])\subset[c,\omega_1[$. Moreover,
  the set $F=\{\gamma\in\omega_1\,:\, \max\{\gamma,d(\gamma)\}=\gamma\}$ is club.
\end{lemma}
\ \\
\parbox{0.6\textwidth}{\noindent
Thus, if $\gamma\in F$, the dark region of the figure on the right is mapped
inside $[\gamma,\omega_1[$. The proof of the existence of $d(c)$ is easy: For
all $b$, $f\restrict{\RR\times\{b\}}$ is cofinal, so by Lemma \ref{lemme1} there is some
$d(b,c)$ such that $f([d(b,c),\omega_1[\times\{b\})\subset[c,\omega_1[$. Take then
$d(c)=\sup_{m\in\omega}d(b_m,c)$ where $\{b_m\}_{m\in\omega}$ is a dense subset of $[0,c]$.
To see that $F$ is club, one shows that $d\restrict{\omega_1}$ is continuous,
see \cite{meszigueshom} for the details.}
\parbox{0.05\textwidth}{\ }
\parbox{0.35\textwidth}{\epsfig{figure=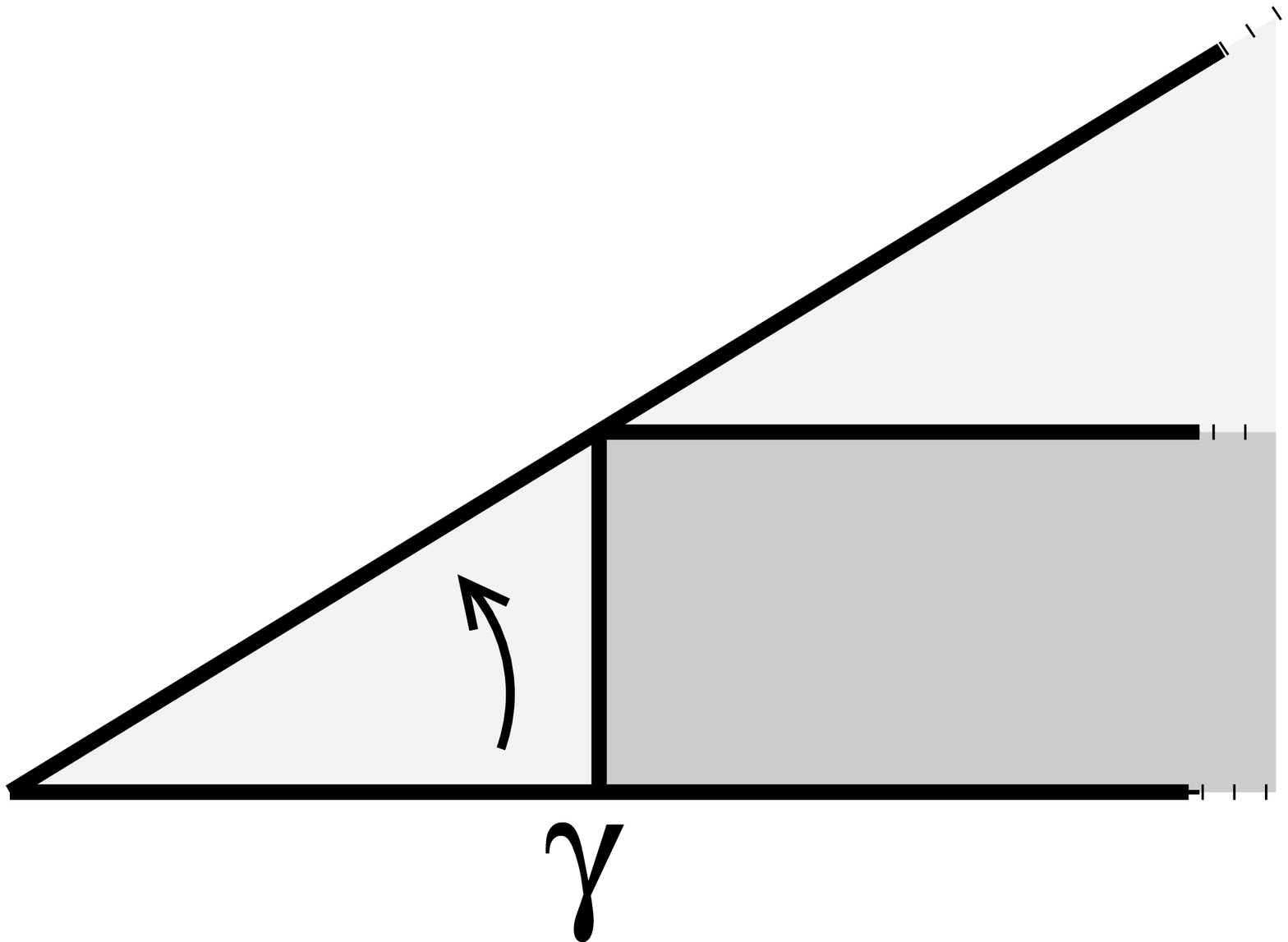,height=2.5cm}}

\begin{lemma}\label{blocks preserving}
  Let $*\in\{d,h\}$, $f\co \mathsf{C}\to\RR$ be continuous and $*$-cofinal. Then,
  $$ 
  E_*=\{\gamma\in\omega_1\,:\, f(A^-_*(\gamma)\subset[0,\gamma]
  \text{ and } f(A^+_*(\gamma)\subset[\gamma,\omega_1[ \} 
  $$
  is club.
\end{lemma}

\proof[Proof of Lemma \ref{blocks preserving}] 
We prove that $E_*^-=\{\gamma\in\omega_1\,:\, f(A^-_*(\gamma)\subset[0,\gamma]\}$ and
$E_*^+=\{\gamma\in\omega_1\,:\,f(A^+_*(\gamma)\subset[\gamma,\omega_1[ \}$ are both club.
Closeness is obvious. The proof that $E_*^-$ is unbounded is easy, we leave it to the reader.
Let thus $\gamma_0\in\omega_1$ be fixed, we show that there is some $\gamma\ge\gamma_0$ in $E_*^+$.
As in Lemma \ref{lemme1}, it is enough to show that for all $\gamma$, there is a $\beta(\gamma)\ge\gamma$ with
$f(A^+_*(\beta(\gamma))\subset[\gamma,\omega_1[$: taking the limit of the sequence
$\gamma_{m+1}=\beta(\gamma_m)$ gives a point in $E_*^+$. So,
suppose absurdly that 
\begin{equation}
  \label{eq1}
  \exists\gamma\,\forall\beta\ge\gamma\,\exists x_\beta\in A_*^+(\beta) 
  \text{ with } f(x_\beta)\le\gamma.
\end{equation}
\ \\
\parbox{0.6\textwidth}{Suppose 
first that $*=d$. 
By (\ref{eq1}) for all $\beta\ge\gamma_0$ there must be a $x_\beta$ in $A_d^+(\beta)$,
the dark region of the 
figure on the right (top),
with $f(x_\beta)\le\gamma$.
Then, starting with $\beta_0=\gamma_0$, we may define sequences $\beta_m$ and $x_{\beta_m}$ as
on the figure on the right (bottom), with $f(x_{\beta_m})\le\gamma$. The sequence 
$x_{\beta_m}$ converges to some point $x=(\beta,\beta)$ on the diagonal $\Delta_d$,
with $\beta\ge\gamma_0$, for which
$f(x)\le\gamma$. Since $\gamma_0$ was arbitrary, it follows that 
$\left(f\restrict{\Delta_d}\right)^{-1}([0,\gamma])$ is unbounded. By Lemma \ref{lemme1}, 
this implies that $f\restrict{\Delta_d}$ is bounded, which is impossible since
$f$ is $d$-cofinal.
}
\parbox{0.05\textwidth}{\ }
\parbox{0.35\textwidth}{\epsfig{figure=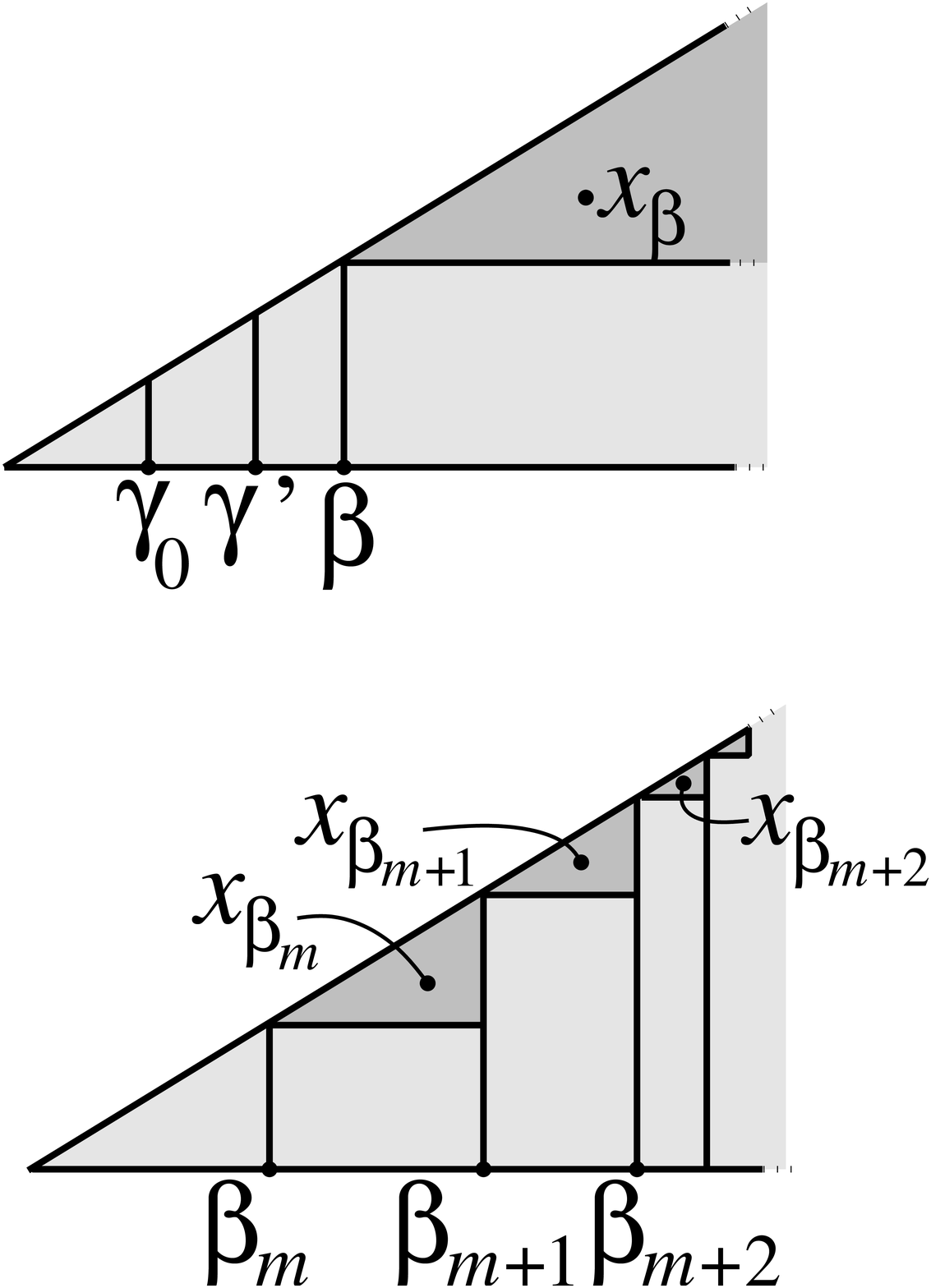,height=5cm}}
\vskip .1cm
If $*=h$, 
Since $f$ is $h$-cofinal, by Lemma \ref{equiv lemme 5.4} there is some
$\gamma'\ge\gamma_0$ with \linebreak
$f([\gamma',\omega_1[\times[0,\gamma'])\subset[\gamma',\omega_1[$, so
there must again be a $x_\beta\in A^+_d(\gamma')$
with $f(x_\beta)\le\gamma$, and we finish as before (using Lemma \ref{h d cof}).
\endproof

\begin{lemma}\label{bounded}
  If $f\co \mathsf{C}\to\RR$ is continuous and $h$-bounded, 
  $$ 
  G=\{\gamma\in\omega_1\,:\, f\restrict{[\gamma,\omega_1[\times\{b\}}
  \text{ is constant }\forall b\in[0,\gamma]\}
  $$
  is club.
\end{lemma}
If $\gamma\in G$, then $f$ is constant on any horizontal at height $\le\gamma$, which are the thin lines of
Figure 2.
The proof of this lemma is very similar to that of Lemma \ref{equiv lemme 5.4}: We
use a dense subset $\{b_m\}_{m\in\omega}$ of $[0,c]$ and Lemma \ref{lemme1} to obtain a
$\wt{d}(c)$ for which $f\restrict{[\wt{d}(c),\omega_1[\times\{b\}}$ is constant for all $b\in[0,c]$.
We omit the details.
\begin{lemma}\label{d bounded}
  If $f\co \mathsf{C}\to\RR$ is continuous and bounded, there is some $\gamma$ such that
  $f\restrict{[\gamma,\omega_1[^2\cap\mathsf{C}}$ 
  and $f\restrict{[\gamma,\omega_1[\times\{b\}}$ are constant $\forall b\le\gamma$.
\end{lemma}
\proof
  Since $f\restrict{\Delta_d}$ is bounded, by Lemma \ref{lemme1} there is some $c$ with
  $f\restrict{\Delta_d\cap[c,\omega_1[^2}\equiv a$ for some $a$. Choose sequences $a_m<a<a_m'$ 
  converging to $a$ and use Lemma \ref{open sets} to conclude that
  $f^{-1}(a)=\bigcap_{m\in\omega}f^{-1}\bigl(]a_m,a_m'[\bigr)\supset]\gamma,\omega_1[^2\,\cap\,\mathsf{C}$ for
  some $\gamma$. Then, apply Lemma \ref{bounded} and increase $\gamma$ if needed to obtain the
  claimed properties for $f\restrict{[\gamma,\omega_1[\times\{b\}}$.
\endproof
We can now prove the partition properties of Figure 2.
If $f$ is $h$-cofinal, Lemma \ref{blocks preserving} enables us to define a strictly increasing sequence
$\left<\beta_\gamma\,:\,\gamma\in\omega_1\right>$ with $\beta_0=0$ and
$\beta_\gamma=\sup_{\gamma' <\gamma}\beta_{\gamma'}$ if $\gamma$
is limit, such that $\beta_\gamma\in E_h$ for all $\gamma$. 
Letting 
\begin{equation}\label{Ph}
P_h(\gamma)=A^+_h(\beta_\gamma)\cap A^-_h(\beta_{\gamma+1}), 
\end{equation}
we get
$f(P_h(\gamma))=[\beta_\gamma,\beta_{\gamma+1}]$ if $\gamma>0$ and
$f(P_h(0))\subset[0,\beta_1]$. If $f$ is $h$-bounded and $d$-cofinal, we may find a similar
$\omega_1$-sequence in $E_d\cap G$, for which
\begin{equation}\label{Pd}
P_d(\gamma)=A^+_d(\beta_\gamma)\cap A^-_d(\beta_{\gamma+1})
\end{equation}
satifies $f(P_d(\gamma))=[\beta_\gamma,\beta_{\gamma+1}]$ for $\gamma>0$,
$f(P_d(0))\subset[0,\beta_1]$,
and $f$ is constant on $[\beta_\gamma,\omega_1[\times\{b\}$
for each $b\in[0,\beta_\gamma[$.
If $f$ is bounded, Lemma \ref{d bounded} gives a $\beta_1$ for which $f$ is constant on the horizontals
$[\beta_1,\omega_1[\times\{b\}$ for $b\ge\beta_1$ and on $[\beta_1,\omega_1[^2\cap\mathsf{C}$.


\section{Proof of Theorem \ref{theorem1}}\label{sec4}
Suppose that $f\co\mathsf{C}\to\RR$ is $h$-cofinal, and let $P_h(\gamma)$ be as in the previous section.
Let $p\co\mathsf{C}\to\RR$ be the projection on the first coordinate. We show that $f$ and $p$ are homotopic.
First, $f(P_h(0))\subset[0,\beta_1]$, we may contract this interval continuously to $\{\beta_1\}$ 
(without moving $[\beta_1,\omega_1[$) and thus assume that $f\restrict{P_h(0)}=p\restrict{P_h(0)}$.
If $\gamma>0$, $f(P_h(\gamma)=[\beta_\gamma,\beta_{\gamma+1}]$ and 
$f\bigl((\{\beta_\gamma\}\times\RR)\cap\mathsf{C}\bigr)=\{\beta_\gamma\}$. Applying Lemma \ref{homotopy} to
$f\restrict{P_h(\gamma)},p\restrict{P_h(\gamma)}$ we obtain homotopies $h_t^\gamma$ between $f$ and $p$ in {\em each}
$P_h(\gamma)$. Letting $h_t(x)$ be $h_t^\gamma(x)$ for $x\in P_h(\gamma)$, we obtain a continuous homotopy between
$f$ and $p$ (notice that $f$ and $p$ agree on $P_h(\gamma)\cap P_h(\gamma+1)$, so $h_t^\gamma$ is constant on these sets).
If $f\co\mathsf{C}\to\RR$ is $h$-bounded and $d$-cofinal, it is homotopic to the projection on the second coordinate:
Find homotopies $h_t^\gamma$ in $P_d(\gamma)$ and then, since $f$ is constant on the horizontals depicted in Figure 2
(middle), we might extend these homotopies in the obvious way to all of $\mathsf{C}$.

It is now easy to generalize this proof to the case of $f,g\co M_{\alpha,s}\to\RR$ satisfying $\mathfrak{C}(f)=\mathfrak{C}(g)$. 
The key point is that since $\alpha<\omega_1$, thanks to Lemma \ref{club} there is an $\omega_1$-sequence
$\left<\beta_\gamma\,:\,\gamma\in\omega_1\right>$ that yields the partition properties 
of Figure 2 for {\em each} $f\restrict{C_\delta},g\restrict{C_\delta}\co C_\delta\to\RR$ ($\delta<\alpha$).
We have pictured the situation for $M_{3,s}$, $s=\left<\uparrow\uparrow\downarrow\right>$ in Figure 3,
firstly with $\mathfrak{C}(f)=\{1,2\}$, and secondly with $\mathfrak{C}(f)=\{2\}$.
The reader is encouraged to have these pictures in mind while going through the proof.
We write $A^+_{*,\delta}(\gamma),A^-_{*,\delta}(\gamma),P_{*,\delta}(\gamma)$ 
for the copies of $A^+_*(\gamma),A^-_*(\gamma),P_*(\gamma)$ in $C_\delta$ ($*\in\{d,h\}$).
As before, we contract $[0,\beta_1]$ to $\{\beta_1\}$
and may thus assume that $f$ and $g$ agree
on $\cup_{\delta<\alpha}A^-_{h,\delta}(\beta_1)$. Looking at Figure 2 ($b$-type), we see that 
if $f,g$ were bounded on $C_\delta$, they become constant after this contraction.
Let now $*(\delta)$ be $h$ if $f\restrict{C_\delta}$ is $h$-cofinal, $d$ if
$f\restrict{C_\delta}$ is $h$-bounded and $d$-cofinal and $b$ if $f\restrict{C_\delta}$ is bounded; the 
partition in $C_\delta$ is thus of the $*(\delta)$-type.
If $*(\delta)=b$, $f\restrict{C_\delta}=g\restrict{C_\delta}\equiv\beta_1$.
There cannot be an infinite (consecutive) sequence of $\delta$ with $*(\delta)=d$,
because it yields sequences 
$\delta_m<\delta'_m<\delta_{m+1}$ ($m\in\omega$) with $f$ $\delta_m$-cofinal and $\delta'_m$-bounded, which
is impossible by continuity (the sequences converge to the same $\delta$, and $f$ cannot be both $\delta$-bounded
and $\delta$-cofinal). 
Moreover, the union of consecutive partition blocks of the $h$-type (preceded and/or
followed by a block of the $d$-type) is homeomorphic to $[0,1]^2$.
Thus,
\begin{equation}{\label{P}}
   P(\gamma)=\bigcup_{\delta<\alpha,*(\delta)\not= b} 
   P_{*(\delta),\delta}(\gamma)
\end{equation}
is homeomorphic to a (finite) disjoint union of squares $[0,1]^2$ when $\gamma\ge 1$,
or to $[0,1]\times\mathbb{S}^1$ (if
$*(\delta)=h$ for all $\delta<\alpha$).
Moreover,
$f(P(\gamma))=[\beta_\gamma,\beta_{\gamma+1}]$ if $\gamma\ge 1$,
$f(P(0))=\{\beta_1\}$, so $f(P(\gamma)\cap P(\gamma+1))=\{\beta_{\gamma+1}\}$, and similarly for $g$.
As above, we define homotopies in each $P(\gamma)$ ($\gamma\ge 1$) and
then, extend it in each $C_\delta$ where $*(\delta)=d$ on the `horizontals' outside $P_{*(\delta),\delta}(\gamma)$.
The homotopy does not move the $C_\delta$ with $*(\delta)=b$, it is therefore continuous.
\endproof
\begin{center}
  \epsfig{figure=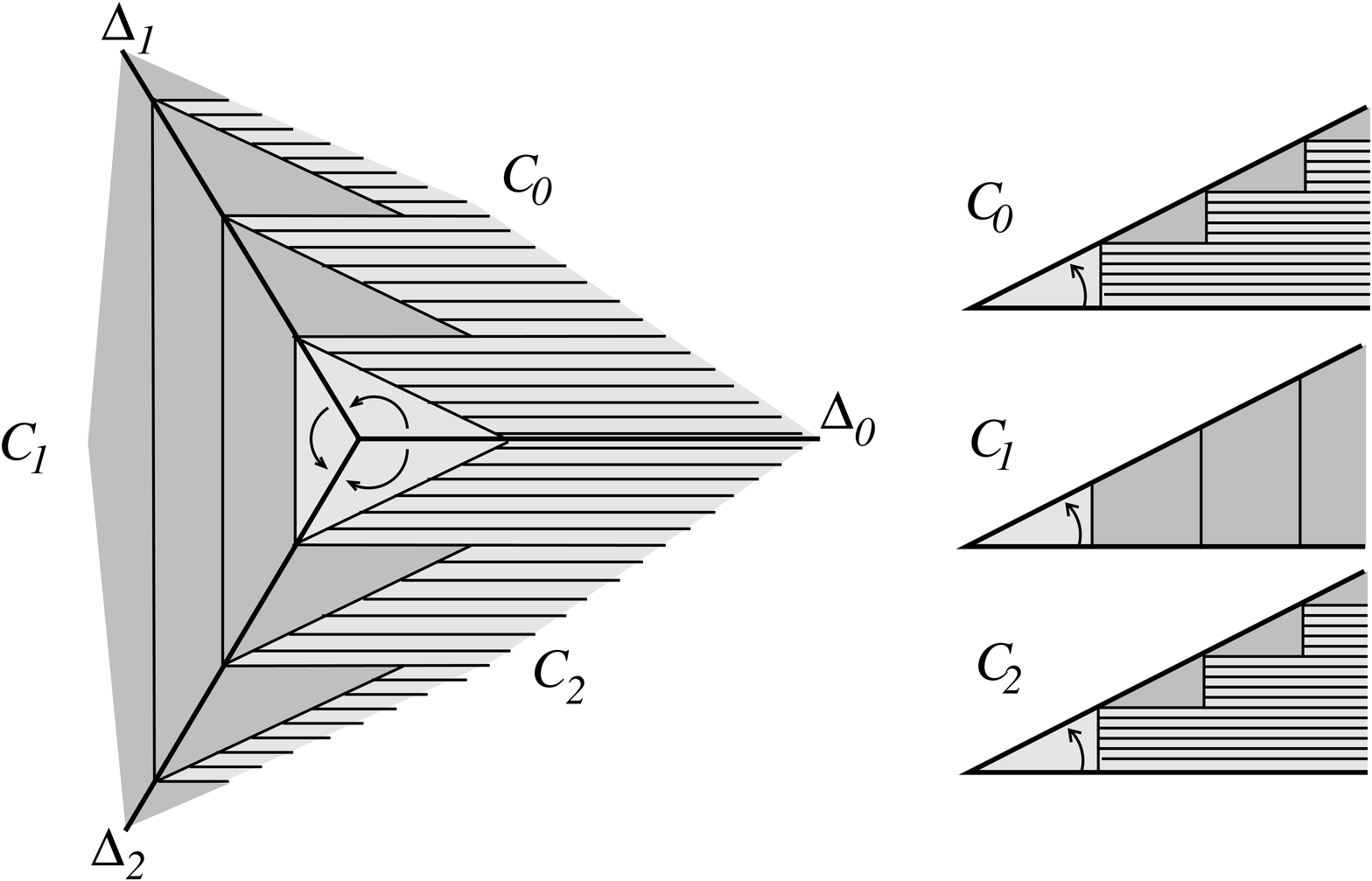,height=6cm}\\
  \epsfig{figure=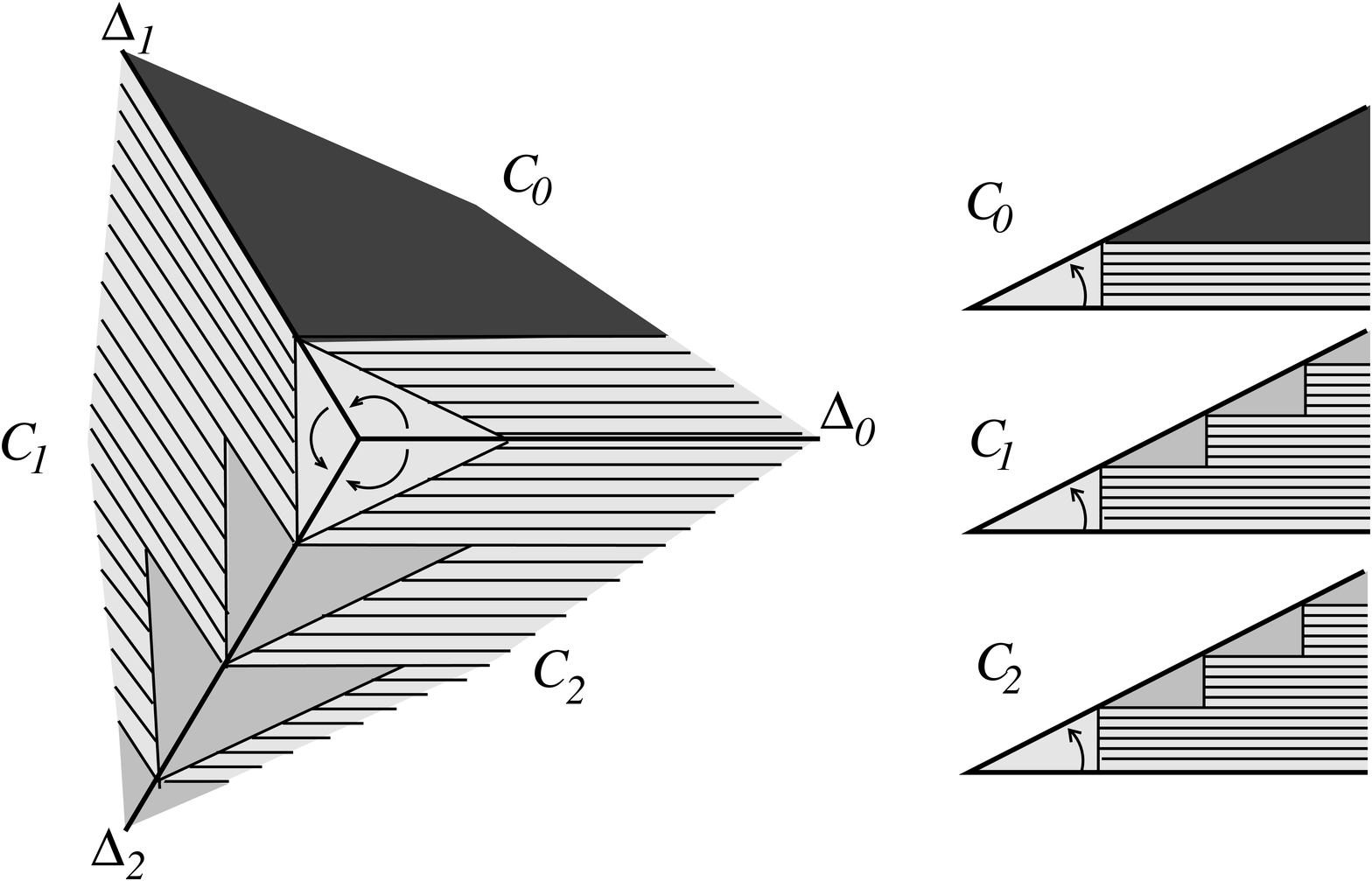,height=6cm}\\
  Figure 3: Partitions for $M_{3,\left<\uparrow\uparrow\downarrow\right>}$, $\mathfrak{C}(f)=\{1,2\},\{2\}$.
\end{center}



\section{Surfaces that do not contain $\RR$}\label{sec6}
Up to now we have built surfaces containing more and more non homotopic copies of $\RR$. We
shall now do the exact opposite: Obtain surfaces $M$ that do not contain $\RR$ but 
which satisfy $[M,\RR]=[\RR,\RR]$. This idea was given to us by A. Henriques, and appears also in
\cite{Nyikos:1992} (concluding remarks).

Since $\RR$ can be given a structure of \CC{\infty} (even analytic) manifold, it has a tangent bundle.
This bundle is nontrivial (see for instance the appendix of \cite{Spivak:vol1}) and depends on the 
smoothing: Nyikos proved in \cite{Nyikos:1992} that there are $2^{\aleph_1}$ non isomorphic
such bundles. 
However, the choice of the smoothing is immaterial here, as we will only use properties 
common to each tangent bundle of $\RR$. \\
\parbox{0.5\textwidth}{So, we denote by
$T$ be the tangent bundle of $\RR$ (given by some smoothing), and by $T_0\subset T$ the $0$ section.
$T\backslash T_0$ falls apart into homeomorphic submanifolds $T^+$ and $T^-$. 
Since $T$ is not trivial, $T^+$ does not contain a copy of $\RR$ (otherwise there is a non $0$ section).
There is a natural action of $\Z$ on (the fibers of) $T$ given by $(x,i)\mapsto 2^{i}\cdot x$.
We define $M$ to be the quotient
$T^+/\Z$, where we moreover identify each point in the fiber above $0$.
(We do this identification to obtain a longplane.)
It is clear that $M$ is an $\omega$-bounded surface.
Away from $0$, we have a covering $\pi\co T^+\to M$ and a (non trivial)
bundle (with fiber $\mathbb{S}^1$) $p\co M\to\RR$. 
We remark that there is no embedding $\RR\to M$: Since $\RR$ is connected, (locally)
pathwise connected with trivial fundamental group, and $\pi$
is a covering, such an embedding could be lifted up into $T^+$ (see \cite{Spanier},
Theorem 5, Chapter 2, Section 4).
Whether $M$ contains $\omega_1$ or not may depend upon
axioms stronger than ZFC, see the concluding remarks in \cite{Nyikos:1992}.}
\parbox{0.45\textwidth}{\ \ \epsfig{figure=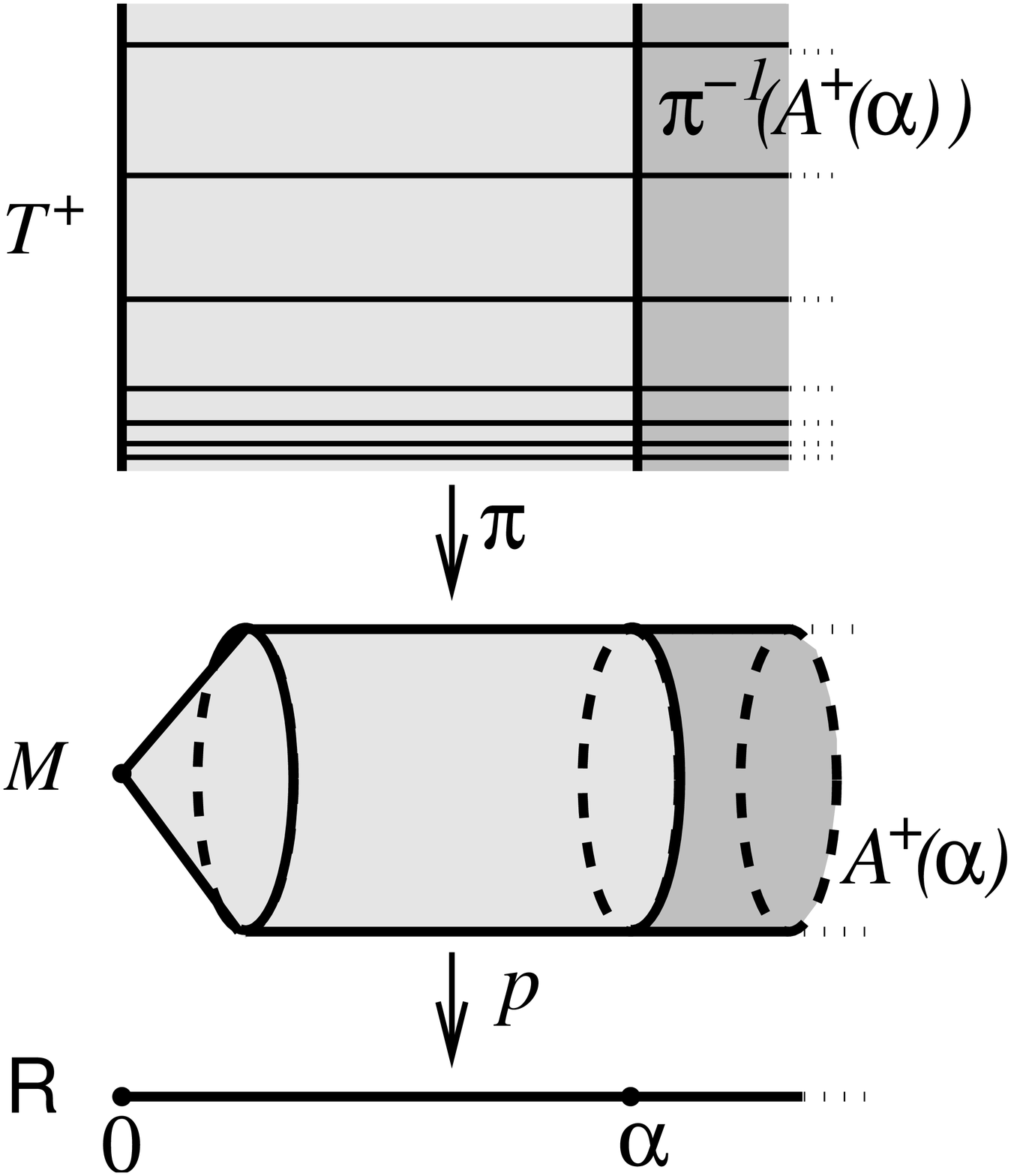,height=6.3cm}}
Let $A^+(\alpha)=p^{-1}([\alpha,\omega_1[)\subset M$, $A^-(\alpha)=p^{-1}([0,\alpha])\subset M$.
We have:
\begin{lemma}\label{partitionT}
If $f\co M\to\RR$ is continuous and unbounded, for each $\beta\in\RR$ there is 
$\alpha(\beta)\in\RR$ with $f(A^+(\alpha(\beta)))\subset [\beta,\omega_1[$.
\end{lemma}
\proof
Suppose otherwise that there is some $\beta$ such that for each $\alpha$, 
$f(A^+(\alpha))\cap [0,\beta]\not=\emptyset$. 
Then $U=f^{-1}([0,\beta +1[)$ and $V=f^{-1}(]\beta +1,\omega_1[)$ are open and
unbounded in $M$, and contain the closed unbounded sets 
$f^{-1}([0,\beta])$ and $f^{-1}([\beta +2,\omega_1[)$. 
If $C\subset M$ is closed and unbounded, it must intersect 
$p^{-1}(\omega_1)$ in a closed unbounded set: Just take a sequence $x_\alpha$ in $C$ such
that $x_\alpha\in p^{-1}([\gamma_\alpha,\gamma_\alpha +1])$ with $\gamma_\alpha>\gamma_{\alpha'}$
whenever $\alpha>\alpha'$; by sequential compacity of $M$, a subsequence converges to some
$x\in p^{-1}(\sup \gamma_\alpha)$.
So, $\pi^{-1}(U)$ and $\pi^{-1}(V)$ are large open sets of $T^+$ (see Definition 4.8 p. 149
of \cite{Nyikos:1992}).
By \cite[Cor. 4.11 p. 151]{Nyikos:1992}, 
$\pi^{-1}(U)\cap\pi^{-1}(V)$ is also large, and thus $U\cap V$ is nonempty, which 
contradicts their definition.
\endproof
\begin{cor}\label{partitionTcor}
  If $f\co M\to\RR$ is continuous and unbounded, then
  $$H=\{\alpha\,:\,f(A^+(\alpha))\subset[\alpha,\omega_1[,\,f(A_\alpha^-)\subset[0,\alpha]\}$$
  is closed unbounded.
\end{cor}
There are therefore at least two homotopy classes of maps $M\to\RR$.
Using homotopies in $P(\alpha)=A^+(\beta_\alpha)\cap A^-(\beta_{\alpha+1})$
(where $\beta_\alpha$ is a strictly increasing closed sequence in $H$), we see as in Section \ref{sec4}
that two unbounded maps
$M\to\RR$ are homotopic, and we get:
\begin{thm}
 $[M,\RR]\simeq\{0,1\}$.
\end{thm}
Notice that if $f\co M\to M$ is continuous and unbounded, Corollary \ref{partitionTcor}
applied to $p\circ f$
yields
$f(P(\alpha))\subset P(\alpha)$ for all $\alpha>0$. 
The entire discussion above could be made for the quotient $M'=T^+/\Z$ (without identifying the points above $0$).
Since $P(\alpha)$ is homeomorphic to
$\mathbb{S}^1\times[0,1]$ whenever $\alpha>0$, the following is plausible:
\begin{prob} Do we have $[M,M]\simeq\Z\cup\{*\}$, where $\{*\}$ is the class of bounded maps,
 and $[M',M']\simeq\Z\times\{0,1\}$ ?
 Does it depend on the smoothing?
\end{prob}
(Using the techniques developped in the previous sections, one can easily prove $[\mathbb{S}^1\times\RR,\mathbb{S}^1\times\RR]\simeq\Z\times\{0,1\}$.) 
A related question (asked by D. Gauld in \cite{Gauld}) is the following. 
Recall that the mapping class group of a topological space $X$ is the group of homeomorphisms of
$X$ up to {\em isotopy}.
\begin{prob}
  Using the partition properties, can we describe explicitly the mapping class group of
  the manifolds described in this paper? Can we do it for any longplane?
\end{prob} 
A third question asks about {\em directions}.
Let $X,Y$ 
be $\omega$-bounded manifolds with canonical 
sequences $\left<U_\alpha\,:\,\alpha\in\omega_1\right>$ and 
$\left<V_\alpha\,:\,\alpha\in\omega_1\right>$ respectively.
We say that a club subset $D\subset X$ is 
a $Y$-direction (in $X$) if
whenever $f\co X\to Y$ is continuous with $f\restrict{D}$ unbounded,
then for every $\alpha\in\omega_1$ there is $\beta(\alpha)\in\omega_1$  
with
$f(D\backslash U_{\beta(\alpha)})\subset Y\backslash V_\alpha$. 
Any copy of $\RR$ (or of $\omega_1$) is an $\RR$-direction by Lemma
\ref{lemme1}, basically the same proof shows that it is in fact a $Y$-direction for all $\omega$-bounded
manifolds $Y$. The surface $M$ of this section is an $\RR$-direction and an $M$-direction, while
$\RR^2$ is {\em not} an $\RR$-direction (take $f(x,y)=\min\{x,y\}$).
\begin{prob}\label{prob:directions1}
  Is there a non-metrizable $\omega$-bounded manifold (surface, longplane) with no $\RR$-direction? 
\end{prob}
\begin{prob} \label{prob:directions2}
  Suppose that $X,Y$ are longplanes (thus with trivial fundamental groups).
  Is $[X,Y]$ described (in any sense) by the $Y$-directions of $X$ (at least when $Y=\RR$)?
\end{prob}
(See also Problem \ref{problem1} below.)
A positive answer to Problem \ref{prob:directions1} would require some axioms beyond ZFC, 
since the proper forcing
axiom implies that any non-metrizable $\omega$-bounded manifold contains a copy of $\omega_1$.
An example not containing $\omega_1$ (using $\diamondsuit$) is given in \cite{Nyikos:1984} (Example 6.9). 
This example
has the property that any closed non compact subset of $X$ contains the bones of the skeleton on a club subset.
(The skeleton is the canonical sequence, the bones are the boundaries $\wb{U}_\alpha\backslash U_\alpha$.)
The reader is invited to check that $X$ is therefore an $\RR$-direction, and 
$[X,\RR]\simeq\{0,1\}$. More generally, using the techniques of the preceding sections one
easily shows the following:

\begin{prop}
   Let $X$ be an $\omega$-bounded manifold which is an $\RR$-direction. Then $[X,\RR]\simeq\{0,1\}$.
\end{prop}


\renewcommand\thesection{\Alph{section}}
\setcounter{section}{0}

\section{Appendix: More $\omega$-bounded surfaces.}\label{sec5}
In this appendix we show that the constructions of Section \ref{sec2} can be generalized to obtain surfaces 
that contain uncountably many mutually non homotopic copies of $\RR$, and obtain results similar to
Theorem \ref{theorem1}. We shall not provide the proofs in detail, since they differ only slightly
from those before.
\esp
We begin by defining the $\omega_1$-octant $\mathsf{C}_{\omega_1,s}$,
for $s\co \omega_1\to\{\uparrow,\downarrow\}$.
We start with $\RR^2$ which we consider as the union of $K_\alpha=\RR\times[\alpha,\alpha+1]$,
$\alpha\in\omega_1$, glued along their boundary components. Then, we replace each $K_\alpha$ by a copy
$\wt{C}_\alpha$ of $\mathsf{C}$ in the way given by $s(\alpha)$ as in Section \ref{sec2} (see Figure 4 below).
Let $I_\alpha$ be a line segment (in $\wt{C}_\alpha$) that joins  
$(\alpha,0)$ to $(\alpha+1,\alpha+1)$ if $s(\alpha)=\uparrow$ and 
$(\alpha,\alpha)$ to $(\alpha+1,0)$ if $s(\alpha)=\downarrow$.
These $I_\alpha$ form 
a `diagonal' $\Delta_{\omega_1}$ in the new space, similar
to $\Delta_d\subset\RR^2$. 
We then define $\mathsf{C}_{\omega_1,s}$ as the surface `below' this diagonal $\Delta_{\omega_1}$, the darker
region in Figure 4, and let $C_\alpha$ be $\mathsf{C}_{\omega_1,s}\cap \wt{C}_\alpha$.
 \begin{center}
  \epsfig{figure=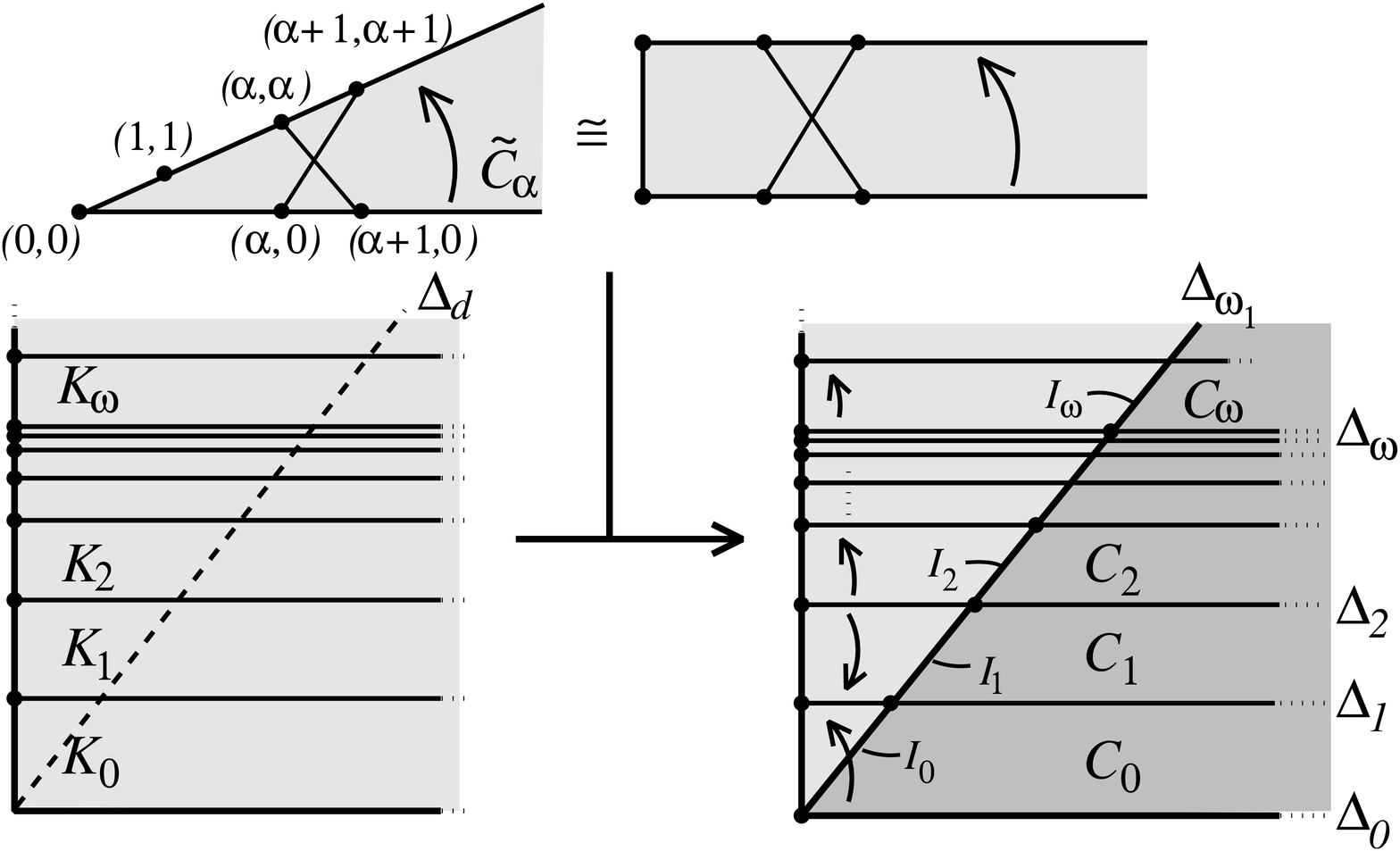,height=5.5cm}\\
  {\em Figure 4: The $\omega_1$-octant $\mathsf{C}_{\omega_1,s}$.}
\end{center}
This time, $\Delta_\alpha$ runs along $\Delta_{\omega_1}$ until it separates from it in $C_\alpha$ (see Figure 4).
$\mathsf{C}_{\omega_1,s}$ shares many properties with
$\mathsf{C}$. In particular, the reader is invited to prove the following:
\begin{lemma}\label{DeltaComega1}
  If $\Delta\subset\mathsf{C}_{\omega_1,s}$ is a copy of $\RR$ that intersects $\Delta_\alpha$
  for all $\alpha<\omega_1$, then $\Delta\cap\Delta_{\omega_1}$ is club. If there is an $\alpha$ such that
  $\Delta\cap\Delta_\gamma=\emptyset$ for all $\gamma\ge\alpha$, then, outside of a compact set,
  $\Delta$ is entirely contained
  in $C_\beta\cup C_{\beta+1}$ for some $\beta\le\alpha$.
\end{lemma}
(Notice that if $\Delta\subset C_\beta$, Lemma \ref{DeltaC} applies.)
We now define $M_{\omega_1,s}$ as $\mathsf{C}_{\omega_1,s}$ where $\Delta_{\omega_1}$ and $\Delta_0$ are 
identified, and
the poset $\mathbb{P}=\mathbb{P}_{\omega_1,s}=\left<\omega_1,\prec\right>$ as in Definition \ref{po}. 
\begin{defi}\label{adapted2}
  $W\subset\mathbb{P}$ is adapted if for all $\alpha<\omega_1$ (strict),
  \begin{itemize}
    \item[i)] $\gamma\in W$, $\gamma\prec\gamma'$ $\Rightarrow$ $\gamma'\in W$,
    \item[ii)] $\forall\beta<\omega_1$ with $\beta$ limit,
    $\exists\gamma(\beta)<\beta$ such that
    $$
      \beta\in W\quad\Leftrightarrow\quad
      \gamma(\beta)\in W \quad\Leftrightarrow\quad \gamma'\in W\,\forall\gamma'\in [\gamma(\beta),\beta[,
    $$
    \item[iii)] If $\exists\beta<\omega_1$ with  
    $\gamma\in W\,\forall\gamma\in [\beta,\omega_1[$, then $0\in W$.
  \end{itemize}
\end{defi}
Compare with Definition \ref{adapted}. Notice that $f\co \mathsf{C}_{\omega_1,s}\to\RR$ could be bounded on each 
$\Delta_\alpha$,
$0<\alpha<\omega_1$ but unbounded on $\Delta_{\omega_1}$ (which is identified with $\Delta_0$ in 
$M_{\omega_1,s}$). This explain why the implication in iii) does not reverse (recall that $\beta\in W$ whenever
$f\restrict{\Delta_\beta}$ is unbounded).

\begin{thmbis}
  $[M_{\omega_1,s},\RR]\simeq\{W\subset\mathbb{P}\,:\,W\text{ adapted}\}$.
\end{thmbis}

\proof[Idea of the proof]
We define $\mathfrak{C}(f)$ for $f\co M_{\omega_1,s}\to\RR$ as before.
First, one checks that there is a bijection between adapted sets and 
$\{\mathfrak{C}(f)\,:\,f\co M_{\omega_1,s}\to\RR\text{ continuous}\}$. Then one shows that $f$ and $g$ are homotopic 
whenever $\mathfrak{C}(f)=\mathfrak{C}(g)$. In the countable case, we used Lemma \ref{club} to show that
there is one club subset of $\omega_1$ for which the various partition properties of Figure 2 hold. Since 
we have to deal with $\aleph_1$ copies of $\mathsf{C}$, this is no longer possible. However, we can proceed as follows.
First, extend $f$ to $\cup_{\delta<\omega_1}\wt{C}_\delta$ in any continuous manner.
Let $E_\delta$ be the club subset of $\omega_1$ such that the desired partition properties (those of
Lemmas \ref{equiv lemme 5.4} to \ref{d bounded}) hold for $f\restrict{\wt{C}_\delta}$, $\beta\in E_\delta$.
Then apply the following classical lemma (see \cite{Kunen}):
\begin{lemma*}
  If $E_\delta\subset\omega_1$, $\delta\in\omega_1$ are club, so is their diagonal intersection
  $$ D=\{\gamma\,:\,\gamma\in E_\delta \, \forall \delta <\gamma\}. $$
\end{lemma*}
We fix an $\omega_1$-sequence $\beta_\gamma$ ($\gamma\in\omega_1$) in $D$ as in Section \ref{sec3}, and define
$P_{*,\delta(\gamma)}$ similarly.  
The reader should be convinced that $C_\delta$ is `not concerned'
by the $P_{*,\delta}$ for $\gamma<\delta$ ($C_\delta$ is `under the diagonal'). 
Letting 
$$ P(\gamma)=\bigcup_{\delta<\gamma, *(\delta)\not= b} P_{*(\delta),\delta}(\gamma) $$
(compare with (\ref{P})), 
we then proceed
as in Section \ref{sec4} to finish the proof.
\endproof

Of course, one could build more and more complicated surfaces, using $\mathsf{C}_{\alpha,s}$ as a building
brick instead of $\mathsf{C}$, and so on, similar results could probably be obtained; however we are not sure it 
is worth the effort. One could also imagine a longplane containing continuously many non homotopic copies
of $\RR$ (such an example was suggested to us by A. Henriques).
The following problems seem more interesting. The first one is a variation on 
Problem \ref{prob:directions2}.

\begin{prob}\label{problem1}
  Let $M$ be a longplane (with trivial homotopy groups), spanned by longlines,
  i.e. for each club (thus non-metrizable) $E\subset M$, there
  is a copy of $\RR$ in $M$ that meets $E$ on a club set.
  Let $\mathcal{Q}$ be the collection of copies of $\RR$ in $M$ up to homotopy. Is there a 
  natural partial order on $\mathcal{Q}$ and a notion of adaptedness such that 
  $[M,\RR]\simeq\{W\subset\mathcal{Q}\text{ adapted}\}$? 
\end{prob} 
We say that a topological space is locally sub-euclidean if each point has a neighborhood that embeds in
$\R^n$.
\begin{prob}
  What conditions are sufficient for a poset $\mathbb{P}$ to find locally sub-euclidean Hausdorff spaces $X,Y$
  with $[X,Y]\simeq \{W\subset\mathbb{P}\text{ adapted}\}$ (for an appropriate notion of adaptedness)?  
\end{prob}
An an exercise, the reader might prove that $|\mathbb{P}|=n<\omega$ is a sufficient condition (built $X$
with $n-1$ copies of $\mathsf{C}$ and take $Y=\RR$). The case $X=Y$ would be interesting.


\esp
{\bf Acknowledgments:} For useful conversations, I wish to thank
  A. Henriques, D. Cimasoni, A. Gabard, D. Gauld, P. de la Harpe and A. Haefliger.

\bibliographystyle{plain}
\bibliography{../biblio}

{\footnotesize
\vskip1cm
\noindent
Mathieu Baillif \\
Section de Math\'ematiques\\
2-4 rue du Li\`evre\\
1211 Gen\`eve 24\\
Switzerland\\
{\tt Mathieu.Baillif@math.unige.ch}}

\end{document}